\documentclass[12pt,letterpaper]{amsart}
\usepackage{amssymb}
\usepackage[dvips]{color}
\usepackage{amsmath}
\usepackage{amsthm}
\usepackage{bbm}
\allowdisplaybreaks[4]

 \newtheorem{thm}{Theorem}[section]
 \newtheorem{cor}[thm]{Corollary}
 \newtheorem{lem}[thm]{Lemma}
 \newtheorem{prop}[thm]{Proposition}
 \theoremstyle{definition}
 \newtheorem{defn}[thm]{Definition}
 \theoremstyle{remark}
 \newtheorem{rem}[thm]{Remark}
 \numberwithin{equation}{section}

 \newcommand{\h}{\mathcal{H}}
 
 \newcommand{\A}{\mathcal{A}}

 \newcommand{\C}{\mathbb{C}}
 \newcommand{\B}{\mathbf{B}}

 \newcommand{\BH}{\mathbf{B}(\mathcal{H})}

 \newcommand{\abs}[1]{\left\vert#1\right\vert}

 \newcommand{\norm}[1]{\left\Vert#1\right\Vert}

\begin{document}
\pagestyle{myheadings}
\title
 {Masas and Bimodule Decompositions of $\rm{II}_{1}$ factors}

\author{ Kunal Mukherjee }

\address{\hskip-\parindent
 Department of Mathematics \\
 Texas A\&M University \\
 College Station TX 77843--3368, USA}

\email{kunal@neo.tamu.edu}

\begin{abstract}
The \emph{measure-multiplicity-invariant} for masas in $\rm{II}_{1}$
factors was introduced in \cite{MR2261688} to distinguish masas that
have the same \emph{Puk\'{a}nszky invariant}. In this paper we study
the measure class in the \emph{measure-multiplicity-invariant}. This
is equivalent to studying the standard Hilbert space as an
associated bimodule. We characterize the type of any masa depending
on the \emph{left-right-measure} using Baire category methods
$($selection principle of Jankov and von Neumann$)$. We present a
second proof of Chifan's result \cite{MR999997} and a measure
theoretic proof of the equivalence of \emph{weak asymptotic
homomorphism property} $($\emph{WAHP}$)$ and singularity that
appeared in \cite{MR999995}.
\end{abstract}

\thanks{}




\dedicatory{}

\commby{}

\maketitle

\section{Introduction}
This is the first of a series of two papers developed by the author
for his Phd thesis. The moral of this paper is : ``\emph{The
phenomena regularity, semiregularity, singularity, weak asymptotic
homomorphism property $($WAHP$)$ and asymptotic homomorphism
property $($AHP$)$ of masas in finite von Neumann algebras can all
be explained by measure theory}". Throughout the entire paper
$\mathcal{M}$ will always
denote a separable $\rm{II}_{1}$ factor. 
Let $A\subset \mathcal{M}$ be a maximal abelian self-adjoint
subalgebra henceforth abbreviated as a \emph{masa}. It is a theorem
of von Neumann that 
$A$ is isomorphic to $L^{\infty}([0,1],dx)$. So the study of masas
in type ${\rm II}_{1}$ factors is understanding its position $($up
to automorphisms$)$ of the ambient von Neumann algebra. For a masa
$A\subset \mathcal{M}$, Dixmier in \cite{MR0059486} defined the
group of \emph{normalizing unitaries} $($or \emph{normaliser}$)$ of
$A$ to be the set
$$N(A)=\left\{u\in \mathcal{U}(\mathcal{M}): uAu^{*}=A\right\},$$
where $\mathcal{U}(\mathcal{M})$ denotes the unitary group of
$\mathcal{M}$. He called \\
$(i)$ $A$ to be \emph{regular} $($also \emph{Cartan}$)$ if $N(A)^{\prime\prime}=\mathcal{M}$,\\
$(ii)$ $A$ to be \emph{semiregular} if $N(A)^{\prime\prime}$ is a
subfactor of $\mathcal{M}$,\\
$(iii)$ $A$ to be \emph{singular} if $N(A)\subset A$.\\
\noindent He also exhibited the presence of all three kinds of masas
in the hyperfinite $\rm{II}_{1}$ factor.\\ 
\indent Two masas $A,B$ of $\mathcal{M}$ are said to be
\emph{conjugate} if there is an automorphism $\theta$ of
$\mathcal{M}$ such that $\theta(A)=B$. If there is an unitary $u\in
\mathcal{M}$ such that $uAu^{*}=B$ then $A$ and $B$ are called
\emph{unitarily} $($\emph{inner}$)$ \emph{conjugate}.\\
\indent Feldman and Moore in \cite{MR0578656}, \cite{MR0578730}
characterized pairs $A\subset \mathcal{M}$, where $A$ is a Cartan
subalgebra, as those coming from $r$-discrete transitive measured
groupoids with a finite measure space X as base. 
It is a remarkable achievement of Connes, Feldman and Weiss
\cite{MR662736} that any countable amenable measured equivalence
relation is generated by a single transformation of the underlying
space. When translated into the language of operator algebras via
the Feldman-Moore construction, this theorem together with a theorem
of Krieger \cite{MR0415341} says that, \emph{if} $\mathcal{M}$
\emph{is any injective von Neumann algebra then any two Cartan
subalgebras are conjugate by an automorphism of} $\mathcal{M}$.
However it follows from their theorem that, there are uncountably
many equivalence classes of Cartan masas up to unitary conjugacy in
the hyperfinite $\rm{II}_{1}$ factor. See \cite{MR663794} for more
examples. There exist $\rm{II}_{1}$ factors with non conjugate
Cartan masas $($see \cite{MR640947}$)$. These masas were
distinguished with the presence or absence of nontrivial
centralizing sequences. Recently Ozawa and Popa have exhibited
examples of $\rm{II}_{1}$ factors with no or at most one
Cartan masa up to unitary conjugacy $($see \cite{MR999999}$)$.\\
\indent The absence of Cartan masas in $\rm{II}_{1}$ factors was
first due to Voiculescu in \cite{MR1371236}. In fact, it was his
amazing discovery that, for any diffuse abelian algebra $A\subset
L(\mathbb{F}_{n})$, the standard Hilbert space
$l^{2}(\mathbb{F}_{n})$ as a $A$, $A$-bimodule contains a copy of
$L^{2}(A)\otimes L^{2}(A)$. His result was improved by Dykema in
\cite{MR1457745} to rule out the presence of
masas in free group factors with finite multiplicity.\\
\indent Getting back to singular masas, in 1960 Puk\'{a}nszky showed
in \cite{MR0113154} that there are countable non conjugate singular
masas in the hyperfinite ${\rm II}_{1}$ factor by introducing an
algebraic invariant for masas in ${\rm II}_{1}$
factors, today known as the \emph{Puk\'{a}nszky invariant}.\\
\indent In 1983 Popa \cite{MR693226} succeeded in showing that all
separable continuous semifinite von Neumann algebras and all
separable factors of type $\rm{III}_{\lambda}$, $0\leq \lambda<1$
have singular masas. Although they exist, citing explicit examples
is a very hard job. In this direction, Smith and Sinclair in
\cite{MR2163938} have given concrete examples of uncountably many
non conjugate singular masas in the hyperfinite $\rm{II}_{1}$
factor. 
White and Sinclair \cite{MR2302742} have given explicit examples of
a continuous path of non conjugate singular masas $($Tauer masas$)$
in the hyperfinite $\rm{II}_{1}$ factor. All the masas in this path
have the same algebraic invariant of Puk\'{a}nszky. 
Subsequently,
White in \cite{MR999998} proved that, any possible value of the
Puk\'{a}nszky invariant can be realized in the hyperfinite
$\rm{II}_{1}$ factor, and any McDuff factor which contains a masa of
Puk\'{a}nszky invariant $\{1\}$ contains masas
of any arbitrary Puk\'{a}nszky invariant.\\
\indent Singularity is often quite hard to check $($see
\cite{MR1132391}$)$. In order to check if a masa is singular
analytical properties ``\emph{AHP}" and ``\emph{WAHP}" were
discovered in \cite{MR1978447}, \cite{MR1904563}. Subsequently
Smith, Sinclair, White and Wiggins in \cite{MR999995} characterized
pairs $A\subset \mathcal{M}$, where $A$ is a singular masa in a
$\rm{II}_{1}$ factor $\mathcal{M}$ to be precisely those for which
$A$ satisfies ``\emph{WAHP}". All the theories that we have outlined
have a common theme namely, ``\emph{What is the structure of the
standard
Hilbert space as a} $w^{*}$ $A$, $A$-\emph{bimodule}.\\
\indent Although many invariants of masas in $\rm{II}_{1}$ are known
the first successful attempt to distinguish masas with a natural
invariant, which have the same Puk\'{a}nszky invariant was due to
Dykema, Smith and Sinclair in \cite{MR2261688}. We call this the
\emph{measure-multiplicity-invariant}. This invariant has two main
components, a measure class and a multiplicity function. This
invariant is not a new
one and has existed in the literature for quite some time. 
For Cartan masas this invariant has very deep meaning and it is very hard to distinguish Cartan masas with this invariant. 
The term
\emph{multiplicity} in the \emph{measure-multiplicity-invariant} is
actually the \emph{Puk\'{a}nszky} invariant of the masa, making it a
stronger invariant. A slightly different invariant was considered by
Neshveyev and St{\o}rmer in \cite{MR1940356}. \\
\indent Our intention is to study singular masas and distinguish
them. In order to do so, it is necessary to think of singularity
from a different point of view. The theory of Cartan masas and
singular masas have so far been viewed from two different angles.
While Cartan masas fit to the theory of orbit equivalence on one
hand \cite{MR0578730}, singular masas fit to the intertwining
techniques of Popa on the other \cite{MR999995}. 
But we would like to have an unique approach that explains all these
phenomena. This is the primary goal of this paper. In this paper, we
characterize masas by studying the structure of the standard
Hilbert space as its associated bimodule. \\
\indent Our second goal is to investigate that, after such a theory
is outlined whether it is possible to obtain proofs of important
theorems regarding masas that were obtained by a number of
researchers by using different ideas. Many old theorems can indeed
be proved but we will mainly prove Chifan's result on tensor
products \cite{MR999997} and the
equivalence of WAHP and singularity \cite{MR999995}. In fact, it seems that
studying the bimodule is the most natural way to approach these problems
as one can exploit a lot of results from Real Analysis.\\
\indent In order to distinguish singular masas which have the same
multiplicity understanding the measure in the
\emph{measure-multiplicity-invariant} is the most important task. So
we study this invariant thoroughly throughout this article. The
second paper will contain explicit calculations of the invariant and questions related to conjugacy of masas.\\
\indent We have learned latter that Popa and Shylakhtenko in
\cite{MR2215135} has results of similar flavor in this direction.
However our way of approaching is completely different. We think
that what is really involved in understanding the types of masas are
the \emph{ measurable selection principle} of \emph{Jankov and von
Neumann} and some generalized version of Dye's theorem on groupoid
normalisers. This is evident from \cite{MR662736}, \cite{MR0578656}
and \cite{MR0578730}. We present completely measure theoretic proofs
based upon Baire category methods $($selection principle$)$. As an
outcome of our approach many theorems related to structure theory of
masas that were proved by different techniques just follows easily
from our
technique.\\ 
\indent Singular masas are often constructed by considering weakly
or strongly mixing actions of infinite abelian groups on finite von
Neumann algebras. We will show that the definition of WAHP can be
strengthened by considering Haar unitaries and Ces\`{a}ro sums which
exactly resembles the definition of
weakly mixing actions. Weakly mixing actions are characterized by null sets of certain measures. The story for singular masas is also similar. \\
\indent This paper is heavily measure theoretic. Much of the measure
theory tools we require are scattered here and there in the
literature. This article is organized as follows. In Sec. 2 we
present some preliminaries of direct integrals, masas and define the
\emph{measure-multiplicity-invariant}. In Sec. 3 we study
disintegration of measures and masas. Sec. 4 deals with generalized
versions of Dye's theorem. Sec. 5 contains the main result i.e the
characterization theorem and a second proof of Chifan's normaliser
formula. This is a very technical section. Sec. 6 contains results
on calculating certain two-norms and a second proof of the
equivalence of WAHP and singularity. Sec. 4 uses the theory of
$L^{1}$ and $L^{2}$ spaces associated to finite von Neumann algebras
for which we have cited related results in that section without
proofs. Appendix A contains structure theorems of measurable
functions satisfying
condition $(N)$ of Lusin which is used in Sec. 5.\\

\noindent \emph{\textbf{Acknowledgements}}: I thank Ken Dykema, my
advisor, for many helpful discussions. I also thank Roger Smith and
David Kerr for providing me with many helpful ideas. I am grateful
to Roger Smith and Allan Sinclair for making their book, ``Finite
von Neumann Algebras and Masas'' available to me for my reference,
much before it reached the bookstore. I would also like to thank
Stuart White for helpful conversations and for suggesting some of
the results which appear in this paper.

\section{Preliminaries}
The paper relies on the theory of direct integrals. So we have
divided this section into three subsections. In the first subsection
we give some well known results about direct integrals of Hilbert
spaces with respect to an abelian von Neumann algebra. In the second
part we give some preliminaries about masas in $\rm{II}_{1}$ factors
and in the third subsection we will define the
\emph{measure-multiplicity-invariant} of masas in $\rm{II}_{1}$
factors. \\
\noindent \textbf{Notation}: Throughout the entire
article $\mathbb{N}_{\infty}$ will denote the set $\mathbb{N}\cup
\{\infty\}$.

\subsection{Direct Integrals}
$ $\\\\
\indent Let a separable Hilbert space $\h$ be the direct integral of
a $\mu$-measurable field of Hilbert spaces $\{\h_{x}\}_{x\in X}$
over the base space $(X,\mu)$ where $X$ is a $\sigma$-compact space
and $\mu$ is a positive, complete Borel measure.
\begin{defn}\label{decomposable_and_diagonalizable}
An operator $T\in \BH$ is said to be \emph{decomposable} relative to
the decomposition $\h\cong\int_{X}^{\oplus}\h_{x}d\mu(x)$ if there
exists a $\mu$-measurable field of operators $T_{x}\in
\mathbf{B}(\h_{x})$, such that $x\mapsto \norm{T_{x}}\in
L^{\infty}(X,\mu)$ and
$T=\int_{X}^{\oplus}T_{x}d\mu(x)$.\\
If $T_{x}=c(x)I_{\h_{x}}$, where $c(x)\in \C$ for almost all $x$,
then $T$ is said to be \emph{diagonalizable}.
\end{defn}
\indent It is easy to see that the fibres of a \emph{decomposable}
operator are uniquely determined up to an almost sure equivalence.
The collection of \emph{diagonalizable} and \emph{decomposable}
operators both form von Neumann subalgebras of $\BH$, with the later
being the commutant of the former. Whenever there is no danger of
confusion we will use the term measurable instead of
$\mu$-measurable.
\begin{thm}\label{fundamental_theorem_for_direct_integrals}
Let $A\subset\BH$ be a diffuse abelian von Neumann algebra on a
separable Hilbert space $\h$. Then there exists a measure space
$(X,\mu)$, where $X$ is a $\sigma$-compact space, $\mu$ is a
positive, Borel, non-atomic, complete measure on $X$ and a
measurable field of Hilbert spaces $\{\h_{x}\}_{x \in X}$, such that
$\h$ is unitarily equivalent to,
\begin{equation}\label{direct_integral_of_Hilbert_spaces}
\h\cong\int_{X}^{\oplus}\h_{x}d\mu(x)
\end{equation}
and $A$ is $($unitarily equivalent to$)$ the algebra of
diagonalizable operators on $\int_{X}^{\oplus}\h_{x}d\mu(x) $ with
respect to this decomposition.
\end{thm}
The dimension function of the decomposition in Thm.
\ref{fundamental_theorem_for_direct_integrals} is defined as
\begin{equation}
\nonumber m : X\mapsto \mathbb{N}_{\infty} \text{ by, } m(x)=
dim(\h_{x}).
\end{equation}
The dimension function $m$ is $\mu$-measurable. Such results are
known in greater generality. For a measure space $(X,\mu)$ we denote
by $[\mu]$ the equivalence class of measures on $X$ that are
mutually absolutely continuous with respect to $\mu$. This
decomposition in Thm. \ref{fundamental_theorem_for_direct_integrals}
and hence the \emph{multiplicity function} is unique up to measure
equivalence from Thm. 3, 4 of Chapter 6 of \cite{MR641217}.\\
\indent \emph{We will be always working with finite measures}. Since
direct integrals of Hilbert spaces does not change when the measures
are scaled, we will most of the time assume that the measures have
total mass $1$. Details of these facts can be found
in \cite{MR1468230}, \cite{MR591683}.\\

\subsection{Basics on Masas in $\rm{II}_{1}$ factors}
$ $\\

\begin{defn}
Given a type $\rm{I}$ von Neumann algebra $\mathcal{B}$ we shall
write Type($\mathcal{B}$) for the set of all those $n$ $\in$
$\mathbb{N_{\infty}}$ such that $\mathcal{B}$ has a nonzero
component of type $\rm{I}_{n}$.
\end{defn}

Let $\mathcal{M}$ be a separable $\rm{II}_{1}$ factor with the
faithful, normal, tracial state $\tau$. This trace induces the
two-norm $\norm{ x }_{2}=\tau(x^{*}x)^{1/2}$ on $\mathcal{M}$ and we
write $L^{2}(\mathcal{M})$ for the Hilbert space completion of
$\mathcal{M}$ with respect to this norm.  Let $\mathcal{M}$ act on
$L^{2}(\mathcal{M})$ via left multiplication. Let $J$ denote the
anti-unitary conjugation operator on $L^{2}(\mathcal{M})$ obtained
by extending the densely defined map $J(x)= x^{*}$.
Inclusions of von Neumann
algebras
will always be assumed to be unital until further notice.\\
\indent Given a von Neumann subalgebra $\mathcal{N}$ of
$\mathcal{M}$, let $\mathbb{E}_{\mathcal{N}}$ be the unique trace
preserving conditional expectation from $\mathcal{M}$ onto
$\mathcal{N}.$ This conditional expectation is obtained by
restricting the orthogonal projection $e_{\mathcal{N}}$ from
$L^{2}(\mathcal{M})$ onto $L^{2}(\mathcal{N})$ to $\mathcal{M}$.\\ 
\indent Let $A \subset \mathcal{M}$ be a masa. Then the augmented
algebra $\A = (A \cup JAJ)^{\prime\prime}$ is an abelian algebra,
with a type I commutant, the commutant being taken in
$\B(L^{2}(\mathcal{M}))$ and the center of $\A^{\prime}$ is $\A$.
The Jones projection $e_{A}$ onto $L^{2}(A)$ lies in $\A$
\cite{MR999996}. Hence, $\A^{\prime}(1 - e_{A})$ decomposes as,
\begin{equation}
\A^{\prime}(1 - e_{A})= \oplus_{n \in\mathbb{N_{\infty}}}\A^{\prime}P_{n}\\
\end{equation}
\noindent where $P_{n}$  $\in$ $\A$ are orthogonal projections
summing up to $1 - e_{A}$  and $\A^{\prime}P_{n}$ is homogenous
algebra of type $n$ whenever $P_{n} \neq 0$.

\begin{lem}\label{diffuse_lemma}
If $A \subset \mathcal{M}$ be a masa and $B \subseteq \mathcal{M}$
be any subalgebra, then $(A \cup JBJ)^{\prime\prime}$ is diffuse.
\end{lem}
\begin{defn}\label{Puk_invariant}
The \emph{Puk\'{a}nszky invariant} of a masa $A$ in $\rm{II}_{1}$
factor $\mathcal{M}$, denoted by $Puk(A)$ (or $Puk_{\mathcal{M}}(A)$
when the containing factor is ambiguous) is $\{$ $n$ $\in$
$\mathbb{N_{\infty}}$ $:$ $P_{n}\neq 0 \}$ which is precisely
Type($\A^{\prime}(1 - e_{A}))$.
\end{defn}
\begin{defn}\label{normalizing_groupoid}
If $A$ is an abelian von Neumann subalgebra of $\mathcal{M}$, let
$\mathcal{GN}(A)$ or $\mathcal{GN}(A,\mathcal{M})$ be the
\emph{normalising groupoid}, consisting of those partial isometries
$v\in \mathcal{M}$ that satisfy $v^{*}v,vv^{*}\in A$ and
$vAv^{*}=Avv^{*}=vv^{*}A$.
\end{defn}
A theorem of Dye \cite{MR0131516} says that, a partial isometry
$v\in \mathcal{GN}(A)$ if and only if there is an unitary $u\in
N(A)$ and a projection $p\in A$ such that $v=up=(upu^{*})u$. Thus
$\mathcal{GN}(A)^{\prime\prime}=N(A)^{\prime\prime}$. Popa in
\cite{MR815434} connected the \emph{Puk\'{a}nszky invariant} to the
type of a masa showing that if $1\not\in Puk( A)$, then $A$ is
singular and that the \emph{Puk\'{a}nszky invariant} of a Cartan
masa is $\{1\}$.

Singularity is difficult to verify. The following two conditions
were introduced in \cite{MR1904563}, \cite{MR1978447} and
\cite{MR999995} as they imply singularity and are often easier to
verify in explicit situations.

\begin{defn}$($Smith, Sinclair$)$\label{AHP} Let $A$ be a masa in a $\rm{II}_{1}$ factor
$\mathcal{M}$.\\
$(i)$  $A$ is said to have the \emph{asymptotic homomorphism
property} $($AHP$)$ if there exists an unitary $v \in A$ such that
\begin{align}
\nonumber \underset{\abs{n}\rightarrow
\infty}\lim\norm{\mathbb{E}_{A}(xv^{n}y)-\mathbb{E}_{A}(x)v^{n}\mathbb{E}_{A}(y)}_{2}=0
\text{ for all }x,y\in \mathcal{M}.
\end{align}
\noindent $(ii)$ A has the \emph{weak asymptotic homomorphism
property} $($WAHP$)$ if, for each $\epsilon >0$ and each finite
subset $x_{1},\cdots, x_{n}\in \mathcal{M}$ there is an unitary
$u\in A$ such that
\begin{align}
\nonumber
\norm{\mathbb{E}_{A}(x_{i}ux_{j}^{*})-\mathbb{E}_{A}(x_{i})u\mathbb{E}_{A}(x_{j}^{*})}_{2}<
\epsilon \text{ for }1\leq i,j\leq n.
\end{align}
\end{defn}

In \cite{MR999995} it was shown that singularity is equivalent to
WAHP. We will prove in Sec. 6 that WAHP is indeed the most natural
property. The next proposition is well known, we state it for
completeness.

\begin{prop}\label{sakai's_theorem_for_injectivity}
Let $\mathcal{N}\subseteq \BH$ be a von Neumann algebra and let
$x_{i,j}\in \text{ } \mathcal{N}$ and $x_{i,j}^{\prime}\in \text{
}\mathcal{N}^{\prime}$ for
$i,j=1,2,\cdots,n$. Then the following conditions are equivalent:\\
$(i)$ $\sum_{k=1}^{n}x_{i,k}x_{k,j}^{\prime}=0$ for all $1\leq i,j\leq n$.\\
$(ii)$ There exist elements $z_{i,j}\in \mathbf{Z}(\mathcal{N})$,
$i,j=1,2,\cdots,n$ such that for all $i,j$
$$\sum_{k=1}^{n}x_{i,k}z_{k,j}=0, \text{ } \sum_{k=1}^{n}z_{i,k}x_{k,j}^{\prime}=x_{i,j}^{\prime}.$$
\end{prop}

\subsection{Measure-Multiplicity-Invariant}
$ $\\

We consider the \emph{conjugacy invariant} for a masa $A$ in a
$\rm{II}_{1}$ factor $\mathcal{M}$ derived from writing the direct
integral decomposition of its left-right action. More precisely, we
choose a compact Hausdorff space $Y$ such that $C(Y)\subset A$, is a
norm separable unital $C^{*}$ subalgebra and $C(Y)$ is \emph{w.o.t}
dense in $A$. $\tau$ restricted to $C(Y)$ gives rise to a
probability measure $\nu$ on $Y$ so that $A$ is isomorphic to
$L^{\infty}(Y,\overline{\nu})$, with $\overline{\nu}$ a completion
of $\nu$. For simplicity of notation we will use the same symbol
$\nu$ to denote its completion. Now $a\otimes b \mapsto aJb^{*}J$,
$a$, $b$ $\in$ $C(Y)$ extends to an injective $*$-homomorphism $\pi$
of $C(Y)\otimes C(Y)$ in $L^{2}(\mathcal{M})$. Indeed, as
$\mathcal{M}$ is a factor so the map,
$$\sum_{i=1}^{n}a_{i}\otimes b_{i} \mapsto
\sum_{i=1}^{n}a_{i}Jb_{i}^{*}J$$ is injective by Prop.
\ref{sakai's_theorem_for_injectivity}. Hence it induces a norm on
$C(Y)\otimes_{alg} C(Y)$. Since abelian $C^{*}$ algebras are
nuclear, this norm must be the min norm, and therefore $a\otimes b
\mapsto aJb^{*}J$ extends to an injective representation of
$C(Y)\otimes C(Y)$ in $L^{2}(\mathcal{M})$. Therefore
$C(Y$$\times$$Y)$ is a \emph{w.o.t} dense unital subalgebra of $\A$,
so that $\A$ is isomorphic to $L^{\infty}(Y \times Y, \eta_{Y\times
Y})$ for a complete, positive, Borel measure $\eta_{Y\times Y}$. By
Lemma \ref{diffuse_lemma}, $\eta_{Y\times Y}$ is non-atomic.
\begin{rem}
In general, if we allow $\mathcal{M}$ to be a finite von Neumann
algebra that is not a factor then the map
$\sum_{i=1}^{n}a_{i}\otimes b_{i} \mapsto
\sum_{i=1}^{n}a_{i}Jb_{i}^{*}J$ is never injective and the measure
will be supported on smaller sets. See Rem. \ref{not_injection}.
This is the reason we consider factors, although most results of
this article goes through even for finite von Neumann algebras. 
\end{rem}
\indent Thus in view of the uniqueness of direct integrals with
respect to an abelian algebra $($see Thm.
\ref{fundamental_theorem_for_direct_integrals}$)$,
$L^{2}(\mathcal{M})$ admits a direct integral decomposition $\{
\h_{x,y}\}$ over the base space ($Y \times Y$, $\eta_{Y\times Y})$
so that $\A\cong L^{\infty}(Y \times Y, \eta_{Y\times Y})$ is the
algebra of \emph{diagonalizable operators} with respect to this
decomposition. Let $m_{Y}$ denote the multiplicity function of the
above decomposition. It is clear from the direct integral
decomposition that, the \emph{Puk\'{a}nszky invariant} of $A$ is the
set of \emph{essential values} of $m_{Y}$ $($also check Cor. 3.2,
\cite{MR1940356}$)$. We will call $[\eta_{Y\times Y}]$ the
\emph{left-right-measure} of $A$. For reasons that will become
clear, we will in most situation use the same terminology for the
class of the measure $\eta_{Y\times Y}$ when restricted to the
\emph{off diagonal}. This will be clear from the context and will
cause no confusion. A related invariant was considered by Neshveyev
and St{\o}rmer in \cite{MR1940356}, which
was a complete invariant for the pair $(A,J)$.\\
\indent Although the existence of such a measure is guaranteed we
need an algorithm to figure out the \emph{left-right-measure}. In
order to do so fix a nonzero vector $\xi$ $\in$
$L^{2}(\mathcal{M})$. The cyclic projection $P_{\xi}$ with range
[$\A\xi$] is in $\A^{\prime}$ and hence \emph{decomposable}. For
$f$, $g$ $\in$ $C(Y)$, there exists a complete positive measure
$\mu_{\xi}$ (we complete it if necessary) on $Y\times Y$ such that
\begin{equation}\label{eta_xi_measure_formulae}
\langle fJg^{*}J\xi, \xi\rangle_{L^{2}(\mathcal{M})}= \int_{Y\times
Y}f(t)g(s)d\mu_{\xi}(t,s).
\end{equation}
$\A$$P_{\xi}$ is a diffuse abelian algebra in
$\B(P_{\xi}(L^{2}(\mathcal{M})))$ with a cyclic vector, so is
maximal abelian. Thanks to von Neumann, we have only one. Therefore,
\begin{equation}\label{eta_xi_measure_decompose}
P_{\xi}(L^{2}(\mathcal{M}))\cong\int_{Y\times
Y}^{\oplus}{\mathbb{C}}_{t,s}d\mu_{\xi}(t,s) \text{ where
}\C_{t,s}=\C.
\end{equation}
\noindent Moreover $\A$$P_{\xi}$ is the \emph{diagonalizable
algebra} with
respect to the decomposition in Eq. \eqref{eta_xi_measure_decompose}.\\
\indent Two orthogonal cyclic subspaces $[\A\xi_{1}]$ and
$[\A\xi_{2}]$ with cyclic vectors $\xi_{1}$, $\xi_{2}$ does not
necessarily keep the fibres of its associated projections
$P_{\xi_{1}}$ and $P_{\xi_{2}}$ orthogonal, neither does assert that
they are direct integrals over disjoint subsets of $Y\times Y$.
However, using the ``gluing lemma" $($Lemma 5.7, \cite{MR2261688}$)$
we single out a measure $\mu_{\xi_{1},\xi_{2}}$ so that
$(P_{\xi_{1}}+P_{\xi_{2}})(L^{2}(\mathcal{M}))$ has a direct
integral decomposition with respect to $(Y\times
Y,\eta_{\xi_{1},\xi_{2}})$ and $\A$$(P_{\xi_{1}}+P_{\xi_{2}})$ is
the \emph{diagonalizable algebra} respecting that decomposition.
This is the step where one will see the possible updates of the
\emph{multiplicity function}. Since we are working on a separable
Hilbert space, after at most a countable infinite iterations of this
procedure we will finally find a measure $\mu$ on $Y\times Y$ so
that
\begin{equation}\label{final_measure_obtained}
L^{2}(\mathcal{M})\cong\int_{Y\times
Y}^{\oplus}\h_{x}^{\prime}d\mu(x)
\end{equation}
\noindent and $\A$ is \emph{diagonalizable} with respect to the
decomposition in Eq. \eqref{final_measure_obtained}. Modulo the
uniqueness of direct integrals 
we have found the
measure. Needless to say, different choices of cyclic subspaces will
produce same measure modulo the uniqueness. However for purpose of
explicit computation to distinguish masas one learns, that nice
choices of cyclic projections $($vectors$)$ is perhaps a
little too costly.\\
\indent For a set $X$ we denote by $\Delta(X)$ the set $\{(x$, $y)$
$\in$ $X$ $\times$ $X$ $: x = y\}$. The restriction of $\tau$ to
$C(Y)\subset A$ gives rise to a Borel probability measure whose
completion is denoted by $\nu_{Y}$.
\begin{lem}\label{properties_of_eta}
The measure $\eta_{Y\times Y}$ has the following properties:\\
$($i$)$ $[\eta_{Y\times Y}]$ is invariant under the flip map
$\theta$ :$(s, t)$ $\mapsto$ $(t, s)$ on $Y \times Y$.\\
$($ii$)$ If $\pi_{1}$ and $\pi_{2}$ denote the coordinate
projections from $Y\times Y$ onto $Y$ then,
\begin{equation}
[({\pi_{i}})_{*}\eta_{Y\times Y} ]= [\nu_{Y}] \text{ for } i =1,2.
\end{equation}
$($iii$)$ The subspace
$\int^{\oplus}_{{\Delta}(Y)}\h_{t,s}d\eta_{Y\times Y}(t,s)$ is
identified with $L^{2}(A)$ and $m_{Y}(t,t)= 1$, $\eta_{Y\times Y}$ a.e. on $\Delta(Y)$.\\
$($iv$)$ The topological $($closed$)$ support of $\eta_{Y\times Y}$ is $Y\times Y$.\\
The multiplicity function $m_{Y}$ has the property that
$$m_{Y}(s,t)=m_{Y}(t,s)$$ almost all $\eta_{Y\times
Y}$.
\end{lem}
Lemma \ref{properties_of_eta} is known so we omit its proof.
Interested readers can consult \cite{MR1940356} or
\cite{MR1000012}. In fact, it is possible to obtain a choice of $\eta_{Y\times Y}$ such that $\eta_{Y\times Y}=\theta_{*}\eta_{Y\times Y}$.\\
\indent  We are now almost ready to give the definition of the
\emph{measure-multiplicity-invariant} of a masa in a separable
$\rm{II}_1$ factor. 
Let $A$ be a masa in $\mathcal{M}$. Let $Y$ be any compact Hausdorff
space such that the unital inclusion of $C(Y)$ in $A$ is
\emph{w.o.t} dense and $C(Y)$ is norm separable. 
To each such Y, we associate a quadruple $(Y,\nu_{Y},[\eta_{Y\times
Y}],m_{Y})$. Define an equivalence relation on the quadruples
$(Y,\nu_{Y},[\eta_{Y\times Y}],m_{Y})$ by\\
$(Y,\nu_{Y},[\eta_{Y\times
Y}],m_{Y})\sim_{m.m}(Y^{\prime},\nu_{Y^{\prime}},[\eta_{Y^{\prime}\times
Y^{\prime}}],m_{Y^{\prime}})$ if and only if there exists a Borel
isomorphism $F :Y \mapsto Y^{\prime}$ such that,
\begin{align}
\nonumber  &F_{*}\nu_{Y} = \nu_{Y^{\prime}},\\
&(F \times F)_{*}[\eta_{Y\times Y}] = [\eta_{Y^{\prime}\times Y^{\prime}}] \text{ and}\\
\nonumber &m_{Y}\circ (F \times F)^{-1} = m_{Y^{\prime}}\text{, }
\eta_{Y^{\prime}\times Y^{\prime}} \text{ a.e.}\nonumber
\end{align}
\noindent We also have, $[\eta_{Y\times Y}]= [\eta_{\mid\Delta(Y)}]$
$+ $ $[\eta_{\mid\Delta(Y)^{c}}]$.\\ Therefore if
$(Y,\nu_{Y},[\eta_{Y\times
Y}],m_{Y})\sim_{m.m}(Y^{\prime},\nu_{Y^{\prime}},[\eta_{Y^{\prime}\times
Y^{\prime}}],m_{Y^{\prime}})$ then,
\begin{align}
(F \times F)_{*}[\eta_{\mid\Delta(Y)^{c}}]&=[\eta_{\mid\Delta(Y^{'})^{c}}],\\
m_{\mid\Delta(Y)^{c}}\circ (F \times F)^{-1}&=
m_{{\mid\Delta(Y^{'})^{c}}},{\text { }} \eta_{\mid\Delta(Y^{'})^{c}}
{\text { a.e}}. \nonumber
\end{align}

\begin{lem}\label{all_classes_equivalent}
If $C(Y_{1}) \subseteq  C(Y_{2})\subset A \subset \mathcal{M}$ be
two w.o.t dense, unital, norm separable $C^{*}$ subalgebras of $A$
then $(Y_{1},\nu_{Y_{1}},[\eta_{Y_{1}\times
Y_{1}}],m_{Y_{1}})\sim_{m.m}(Y_{2},\nu_{Y_{2}},[\eta_{Y_{2}\times
Y_{2}}],m_{Y_{2}})$.
\end{lem}
\begin{proof}
The inclusion $i:C(Y_{1}) \hookrightarrow C(Y_{2})$ results from a
continuous surjection $\theta :Y_{2}\mapsto Y_{1}$. Therefore for
all $f \in$ $C(Y_{1})$,
$$\tau(f)=
\int_{Y_{1}}fd\nu_{Y_{1}}=\int_{Y_{2}}i(f)d\nu_{Y_{2}}=\int_{Y_{2}}(f\circ\theta)d\nu_{Y_{2}}=\int_{Y_{1}}fd(\theta_{*}\nu_{Y_{2}}).$$
Therefore, $\theta_{*}\nu_{Y_{2}}=\nu_{Y_{1}}$.\\
\noindent The inclusion $i$ preserves least upper bounds at the
level of continuous functions. So $i$ extends to a surjective
$*$-homomorphism $\tilde{i}$ between $L^{\infty}(Y_{1},\nu_{Y_{1}})$
and $L^{\infty}(Y_{2},\nu_{Y_{2}})$ which is normal $($Lemma 10.1.10
\cite{MR1468230}$)$. It is easy to see that $\tilde{i}$ is also
implemented by $\theta$. That $\tilde{i}$ is injective is obvious.
So $\theta$ is a Borel isomorphism between the underlying measure
spaces.\\
\noindent Arguing similarly it is easy to see that $\theta\times
\theta :Y_{2}\times Y_{2}\mapsto Y_{1}\times Y_{1}$ implements an
isomorphism between $L^{\infty}(Y_{1}\times Y_{1},\eta_{Y_{1}\times
Y_{1}})$ and $L^{\infty}(Y_{2}\times Y_{2},\eta_{Y_{2}\times
Y_{2}})$. The statements regarding the measure classes now follows
easily.\\
\noindent The statement about the multiplicity function is obvious
from the uniqueness of direct integrals in Thm.
\ref{fundamental_theorem_for_direct_integrals} and the fact
$L^{\infty}(Y_{1}\times Y_{1},\eta_{Y_{1}\times Y_{1}})\cong
L^{\infty}(Y_{2}\times Y_{2},\eta_{Y_{2}\times Y_{2}})\cong \A$.
\end{proof}

\begin{prop}
Let $A\subset\mathcal{M}$ be a masa. The collection of quadruples
$(Y,\nu_{Y},[\eta_{Y\times Y}],$ $m_{Y})$ for $Y$ a compact
Hausdorff space such that $C(Y)\subset A$ is unital, norm separable
and w.o.t dense in $A$, under the equivalence relation $\sim_{m.m}$
has exactly one equivalence class.
\end{prop}
\begin{proof}
If $C(Y_{1})$, $C(Y_{2})\subset A$ be two \emph{w.o.t} dense,
unital, norm separable subalgebras of $A$ then $C^{*}(C(Y_{1})\cup
C(Y_{2}))\cong C(Y_{3})$ for a compact Hausdorff space $Y_{3}$, and
$C(Y_{3})$ is unital, norm separable and \emph{w.o.t} dense in $A$.
Therefore by Lemma \ref{all_classes_equivalent},
$(Y_{3},\nu_{Y_{3}},[\eta_{Y_{3}\times Y_{3}}],m_{Y_{3}})\sim_{m.m}
(Y_{i},\nu_{Y_{i}},[\eta_{Y_{i}\times Y_{i}}],m_{Y_{i}})$ for $i=$1,
2.
\end{proof}

\begin{defn}\label{define_mm}
Let $A\subset \mathcal{M}$ be a masa. We define the
\emph{measure-multiplicity-invariant} of $A$ as the
\emph{equivalence class} of the quadruples $(Y, \nu_{Y},
[\eta_{\mid\Delta(Y)^{c}}],
m_{\mid\Delta(Y)^{c}})$ under $\sim_{m.m}$ where,\\
$(i)$ $Y$ is a compact Hausdorff space such that $C(Y)$ is an
unital,
norm separable and \emph{w.o.t} dense subalgebra of $A$.\\
$(ii)$ $\nu_{Y}$ is the completion of the probability measure
obtained from restricting $\tau$ on $C(Y)$.\\
$(iii)$ $[\eta_{\mid\Delta(Y)^{c}}]$ is the equivalence class of the
measure $\eta_{Y\times Y}$ restricted to $\Delta(Y)^{c}$,\\
$(iv)$ $m_{\mid\Delta(Y)^{c}}$ is the multiplicity function
restricted
to $\Delta(Y)^{c}$,\\
obtained from the \emph{direct integral decomposition} of
$L^2(\mathcal{M})$ over the base space $(Y$ $\times$ $Y,
\eta_{Y\times Y})$ so that $\A$ is the algebra of
\emph{diagonalizable operators} with respect to this decomposition.
\end{defn}

The \emph{measure-multiplicity-invariant} is an \emph{invariant} for
masas in the following sense. If $A \subset \mathcal{M}$ and $B
\subset \mathcal{N}$ are masas in $\rm{II}_1$ factors $\mathcal{M},
\mathcal{N}$ respectively, and there is an unitary $U :
L^{2}(\mathcal{M}) \mapsto L^{2}(\mathcal{N})$ such that, $UAU^{*} =
B$ and $UJ_{\mathcal{M}}AJ_{\mathcal{M}}U^{*} =
J_{\mathcal{N}}BJ_{\mathcal{N}}$ then for any choice of compact
Hausdorff spaces $Y_{A}, Y_{B}$ with\\
$\overline{C(Y_{A})}^{s.o.t}= A$ and $\overline{C(Y_{B})}^{s.o.t}=
B$, $1_{\mathcal{M}}\in C(Y_{A})$, $1_{\mathcal{N}}\in C(Y_{B})$ and
$C(Y_{A}), C(Y_{B})$ norm separable, there exists a Borel
isomorphism
\begin{align}
\nonumber F_{Y_{A},Y_{B}}: (Y_{A},\nu_{Y_{A}}) \mapsto
(Y_{B},\nu_{Y_{B}})\text{ such that,}
\end{align}
\begin{align}\label{explain_mm_invariant}
\nonumber &(F_{Y_{A},Y_{B}})_{*}\nu_{Y_{A}} = \nu_{Y_{B}},\\
&(F_{Y_{A},Y_{B}} \times
F_{Y_{A},Y_{B}})_{*}[\eta_{\mid\Delta(Y_{A})^{c}}]=
[\eta_{\mid\Delta(Y_{B})^{c}}] \text{ and}\\
&m_{\mid\Delta(Y_{A})^{c}}\circ (F_{Y_{A},Y_{B}} \times
F_{Y_{A},Y_{B}})^{-1}= m_{\mid\Delta(Y_{B})^{c}},{\text { }}
\eta_{\mid\Delta(Y_{B})^{c}} {\text { a.e}}. \nonumber
\end{align}
We will denote the \emph{measure-multiplicity-invariant} of a masa
$A$ by $m.m(A)$ (or $m.m_{\mathcal{M}}(A)$ when the containing
factor is ambiguous).\\

\section{Conditional measures and Masas}

\indent As we will see latter, the
\emph{measure-multiplicity-invariant} contains substantial
information of the masa. In order to extract more
information we need to establish some house keeping results in measure theory.\\
\indent Disintegration of measures is a very useful tool in ergodic
theory, in the study of conditional probabilities and descriptive
set theory. \emph{Measurable selection principle} is a term closely
linked with disintegration of measures and has been studied by a
number of mathematicians in the last century. A detailed exposition
of the
existence of disintegration can be found in \cite{MR1484954}.\\
\indent For  the general definition of disintegration of measures we
will restrict to the following set up. Let $T$ be a measurable map
from $(X,\sigma_{X})$ to $(Y,\sigma_{Y})$ where
$\sigma_{X},\sigma_{Y}$ are $\sigma$-algebras of subsets of $X,Y$
respectively. Let $\lambda$ be a $\sigma$-finite measure on
$\sigma_{X}$ and $\mu$ a $\sigma$-finite measure on $\sigma_{Y}$.
Here $\lambda$ is the measure to be disintegrated and $\mu$ is often
the push forward measure $T_{*}\lambda$, although other
possibilities for $\mu$ is allowed.
\begin{defn}\label{definition_of_disintegration}
We say that $\lambda$ has a disintegration $\{\lambda_{t}\}_{t\in
Y}$ with
respect to $T$ and $\mu$ or a $(T,\mu)$ disintegration if:\\
$(i)$ $\lambda_{t}$ is a $\sigma$-finite measure on $\sigma_{X}$
concentrated on $\{T=t\}$ $($or $T^{-1}\{t\})$, i.e.
$\lambda_{t}(\{T\neq t\})=0$, for
$\mu$-almost all $t$,\\
and for each nonnegative measurable function $f$ on $X$\\
$(ii)$ $t\mapsto \lambda_{t}(f)$ is measurable.\\
$(iii)$
$\lambda(f)=\mu^{t}(\lambda_{t}(f))\overset{\text{defn}}=\int_{Y}\lambda_{t}(f)d\mu(t)$.
\end{defn}
\indent In probability theory the measures $\lambda_{t}$ are called
the disintegrating measures and $\mu$ is the mixing measure. One
also writes $\lambda(\cdot\mid T=t)$ for $\lambda_{t}(\cdot)$ on
occasion.\\
\indent When $\lambda$ and almost all $\lambda_{t}$ are probability
measures one refers to the disintegrating measures as $($regular$)$
conditional distributions and $t\mapsto \lambda_{t}$ is called the transition kernel.\\
\indent The reader should be cautious that \emph{``measurable"} in
Defn. \ref{definition_of_disintegration} $(ii),(iii)$ means
measurable with respect to the $\sigma$-algebra of completion of
$\lambda$.
\begin{thm}\cite{MR1484954}$($Existence Theorem$)$\label{Existence_of_disintegration} Let $\lambda$ be a
$\sigma$-finite Radon measure on a metric space $X$ and $T$ be a
measurable map into $(Y,\sigma_{Y})$. Let $\mu$ be a $\sigma$-finite
measure on $\sigma_{Y}$ such that $T_{*}\lambda \ll \mu$. If
$\sigma_{Y}$ is countably generated and contains all singleton sets
$\{t\}$, then $\lambda$ has a $(T,\mu)$ disintegration. The measures
$\lambda_{t}$ are uniquely determined up to an almost sure
equivalence: if $\lambda_{t}^{*}$ is another $(T,\mu)$
disintegration then $\mu(\{t:\lambda_{t}\neq \lambda_{t}^{*}\})=0$.
\end{thm}
The condition $T_{*}\lambda \ll \mu$ in Thm.
\ref{Existence_of_disintegration} is actually necessary for the
disintegration to exist. The original version of Thm.
\ref{Existence_of_disintegration} is due to von Neumann.

\begin{prop}\label{atomic_part_measurable_set}
Let $\lambda$ be a Radon measure on a compact metric space $X$ and
$T$ be a measurable map into $(Y,\sigma_{Y})$. Let $\mu$ be a
$\sigma$-finite measure on $\sigma_{Y}$ such that $T_{*}\lambda \ll
\mu$. Assume that $\sigma_{Y}$ is countably generated and contains
all singleton sets. Let $t\mapsto \lambda_{t}$ denote the $(T,\mu)$
disintegration of $\lambda$. Let $X_{a}$ denote the set of atoms of
$\{\lambda_{t}\}_{t\in Y}$ i.e.
$$X_{a}=\left\{x\in X\mid\text{ } \exists\text{ } t\in Y:
\lambda_{t}(\{x\})>0\right\}.$$ Then $X_{a}$ is a measurable set,
measurable with respect to the $\sigma$-algebra of the completion of
$\lambda$.
\end{prop}
\begin{proof}
There is a measurable set $E\subseteq Y$ with $\mu(E^{c})=0$ such
that for $t\in E$, $\lambda_{t}$ is concentrated on the set
$\{T=t\}$. 
We can assume without loss of generality that $E=Y$.
Now for $t\in Y$, the measure $\lambda_{t}$ is concentrated on
$\{T=t\}$, so\\
$$\left\{x\in X\mid\text{ } \exists\text{ } t\in Y:
\lambda_{t}(\{x\})>0\right\}=\left\{x\in X\mid
\lambda_{Tx}(\{x\})>0\right\}.$$
Let $\mathcal{B}$ be a countable base for the topology on $X$. Then
$$\left\{x\in X\mid
\lambda_{Tx}(\{x\})>0\right\}=\cup_{n=1}^{\infty}X_{a}^{(n)}\text{
with }$$
\begin{align}
\nonumber X_{a}^{(n)}&=\left\{x\in X\mid\text{ }\forall U\in\mathcal{B}: x\in U\Rightarrow \lambda_{Tx}(U)\geq\frac{1}{n} \right\}\\
\nonumber &=\underset{U\in \mathcal{B}}\bigcap\left( {(X\setminus
U)\cup \left\{x\in U\mid \lambda_{Tx}(U)\geq\frac{1}{n}
\right\}}\right).
\end{align}
Therefore, $\left\{x\in X\mid\text{ } \exists\text{ } t\in Y:
\lambda_{t}(\{x\})>0\right\}$ is a measurable set by property $(ii)$
of disintegration. 
\end{proof}
\indent The next few lemmas are undoubtedly known to probablists but
we lack the reference. So we record them for convenience. We will
omit their proofs. For details check \cite{MR1000012}.

\begin{lem}\label{fibre_addition_formulae}
Let $\lambda_{1},\lambda_{2}$ be two Radon measures on a compact
metric space $X$ and $T$ be a measurable map into $(Y,\sigma_{Y})$.
Let $\mu$ be a $\sigma$-finite measure on $\sigma_{Y}$ such that
$T_{*}\lambda_{1},T_{*}\lambda_{2} \ll \mu$. Assume $\sigma_{Y}$ is
countably generated and contains all singleton sets $\{t\}$. Let
$\lambda_{t}^{1},\lambda_{t}^{2}$ be the $(T,\mu)$ disintegration of
$\lambda_{1},\lambda_{2}$ respectively. Let $\lambda_{t}^{0}$ be the
$(T,\mu)$ disintegration of $\lambda_{1}+\lambda_{2}$. Then
$$\lambda_{t}^{0}=\lambda_{t}^{1}+\lambda_{t}^{2}\text{ -}\mu \text{ a.e.}$$
\end{lem}

\begin{lem}\label{fibre_tensor_formulae}
Let $\lambda_{1},\lambda_{2}$ be two Radon measures on compact
metric spaces $X,Y$ and $T,S$ be measurable maps from $X,Y$ into
$(Z,\sigma_{Y})$, $(W,\sigma_{W})$ respectively. Let $\mu,\nu$ be
$\sigma$-finite measures on $\sigma_{Y},\sigma_{W}$ respectively
such that
$T_{*}\lambda_{1}\ll\mu,S_{*}\lambda_{2} \ll \nu$.\\
Assume $\sigma_{Y},\sigma_{W}$ are countably generated and contains
all singleton sets $\{t\},\{s\}$ respectively. Let
$\lambda_{t}^{1},\lambda_{s}^{2}$ be the $(T,\mu),(S,\nu)$
disintegration of $\lambda_{1},\lambda_{2}$ respectively. Let
$\lambda_{t,s}^{0}$ be the $(T\otimes S,\mu\otimes\nu)$
disintegration of $\lambda_{1}\otimes\lambda_{2}$. Then
$$\lambda_{t,s}^{0}=\lambda_{t}^{1}\otimes\lambda_{s}^{2}\text{ -}\mu\otimes\nu \text{ a.e.}$$
\end{lem}

\begin{lem}\label{equivalence_of_measure_imply_equivalence_of_fibre_almost_everywhere}
Let $\lambda_{1},\lambda_{2}$ be two Radon measures on a compact
metric space $X$ and $T$ be a measurable map into $(Y,\sigma_{Y})$.
Let $\mu$ be a $\sigma$-finite measure on $\sigma_{Y}$ such that
$T_{*}\lambda_{1} \ll \mu$ and $T_{*}\lambda_{2} \ll \mu$. Assume
$\sigma_{Y}$is countably generated and contains all singleton sets
$\{t\}$. Let $\lambda_{t}^{1},\lambda_{t}^{2}$ be the
$(T,\mu)$ disintegrations of $\lambda_{1},\lambda_{2}$ respectively.\\
$(i)$ Assume that $\lambda_{1}\ll\lambda_{2}\ll\lambda_{1}$. 
Then for $\mu$ almost all $t$,
$\lambda_{t}^{1}\ll\lambda_{t}^{2}\ll\lambda_{t}^{1}$. Moreover, if
$g=\frac{d\lambda_{1}}{d\lambda_{2}}$ then
$\frac{d\lambda^{1}_{t}}{d\lambda^{2}_{t}}=g_{t}$ a.e. $\mu$, where
\begin{equation}
\nonumber g_{t}=\begin{cases}
                  g_{\mid \{T=t\}} \text{ on }\{T=t\},\\
                  0 \text{  otherwise.}
                \end{cases}
\end{equation}
Conversely if $\lambda_{t}^{1}\ll\lambda_{t}^{2}\ll\lambda_{t}^{1}$
for $\mu$ almost all $t$ 
then
$\lambda_{1}\ll\lambda_{2}\ll\lambda_{1}$.\\
$(ii)$ If $\lambda_{1}\perp\lambda_{2}$ then $\lambda^{1}_{t}\perp
\lambda^{2}_{t}$ for $\mu$ almost all $t$.
\end{lem}

\begin{lem}\label{atom_shift_lemma_under_flip_map}
Let $\lambda$ be a Radon measure on $X\times X$ where $X$ is a
compact metric space. Let $\mu$ be a $\sigma$-finite measure on $X$
such that ${(\pi_{i})}_{*}\lambda\ll\mu$ where $\pi_{i}$, $i=1,2$
are coordinate projections onto $X$.\\
Assume that $\lambda$ is invariant under the flip of coordinates
i.e. $\theta_{*}\lambda\ll\lambda\ll\theta_{*}\lambda$, where
$\theta: X\times X\mapsto X\times X$ by $\theta(x,y)=(y,x)$. Let
$\lambda_{s}^{1},\lambda_{t}^{2}$ be the
$(\pi_{1},\mu),(\pi_{2},\mu)$ disintegrations of $\lambda$
respectively. Then for $\mu$ almost all $t$,
$$\lambda^{1}_{t}\ll\theta_{*}\lambda_{t}^{2}\ll
\lambda^{1}_{t}.$$ In particular, if for $\mu$ almost all $t$,
$\lambda_{t}^{2}$ has an atom at $(s,t)$, then $\lambda_{t}^{1}$ has
an atom at $(t,s)$ almost everywhere.
\end{lem}

\begin{thm}\label{fibre_invariance_theorem}
Let $A\subset \mathcal{M}$ and $B\subset \mathcal{N}$ be masas in
separable $\rm{II}_{1}$ factors $\mathcal{M},\mathcal{N}$. Let
$C(X_{1})\subset A$, $C(X_{2})\subset B$ be w.o.t dense, norm
separable, unital subalgebras of $A,B$ respectively, where $X_{i}$
are compact metric spaces for $i=1,2$. Let $\nu_{X_{i}}$ denote the
tracial measures with respect to the w.o.t dense subalgebras on
$X_{i}$ respectively for $i=1,2$. Let $[\lambda_{1}],[\lambda_{2}]$
denote the left-right-measures of $A$ and $B$ respectively. All the
mentioned measures are assumed to be complete. Suppose there is an
unitary $U:L^{2}(\mathcal{M})\mapsto L^{2}(\mathcal{N})$ such that
$UAU^{*}=B$ and
$UJ_{\mathcal{M}}AJ_{\mathcal{M}}U^{*}$ $=J_{\mathcal{N}}BJ_{\mathcal{N}}$.\\
Then there exists isomorphism of measure spaces $F:X_{1}\mapsto
X_{2}$ such that, $F_{*}\nu_{X_{1}}=\nu_{X_{2}}$ and the following
is
true:\\
Denoting by $\lambda_{t}^{1,X_{1}},\lambda_{s}^{2,X_{1}}$ the
$(\pi_{1},\nu_{X_{1}}),(\pi_{2},\nu_{X_{1}})$ disintegrations of
$\lambda_{1}$ respectively and
$\lambda_{t^{\prime}}^{1,X_{2}},\lambda_{s^{\prime}}^{2,X_{2}}$ the
$(\pi_{1},\nu_{X_{2}})$, $(\pi_{2},\nu_{X_{2}})$ disintegrations of
$\lambda_{2}$ respectively, one has
\begin{align}
\nonumber[\lambda_{t^{\prime}}^{1,X_{2}}]&=[(F\times F)_{*}{\lambda}_{F^{-1}t^{\prime}}^{1,X_{1}}]\text{, }\nu_{X_{2}}\text{ almost all} \text{ }t^{\prime},\\
\nonumber[\lambda_{s^{\prime}}^{2,X_{2}}]&=[(F\times
F)_{*}{\lambda}_{F^{-1}s^{\prime}}^{2,X_{1}}]\text{,
}\nu_{X_{2}}\text{ almost all} \text{ }s^{\prime},
\end{align}
where
$\pi_{1},\pi_{2}$ denotes the projection onto the first and second
coordinates respectively.
\end{thm}

If $(X,\sigma)$ be a measurable space and $\mu$ is a signed measure
on $X$ then we denote by $\norm{\mu}_{t.v}$ to be the total
variation norm of $\mu.$
The next Lemma is used in this
paper but it will be of significant use for computation in the next
paper.

\begin{lem}\label{close_measure_lemma}
Let $\lambda_{n},\lambda,\lambda_{0}$ be Radon measures on a compact
metric space $X$ such that, $\lambda_{0}\neq 0$,
$\lambda_{n}\ll\lambda$ for $n=1,2,\cdots$, $\lambda_{0}\ll\lambda$
and $\lambda_{n}\rightarrow \lambda_{0}$ in $\norm{\cdot}_{t.v}$.
Let $T$ be a measurable map into $(Y,\sigma_{Y})$. Let $\mu$ be a
$\sigma$-finite measure on $\sigma_{Y}$ such that $T_{*}\lambda \ll
\mu$. Assume $\sigma_{Y}$is countably generated and contains all
singleton sets $\{t\}$. Let
$\lambda_{t}^{n},\lambda_{t}^{0},\lambda_{t}$ be the $(T,\mu)$
disintegrations of $\lambda_{n},\lambda_{0},\lambda$ respectively.\\
$(i)$ Then there is a $\mu$ null set $E$ and a subsequence
$\{n_{k}\}$ $(n_{k}<n_{k+1}$ for all $k)$ such that for all $t\in
E^{c}$,
$$\underset{A\subseteq \{T=t\},A\text{ Borel}}\sup
\abs{\lambda_{t}^{n_{k}}(A)-\lambda_{t}^{0}(A)}\rightarrow 0 \text{ as }k\rightarrow \infty.$$\\
$(ii) $Moreover, if for $\mu$ almost all $t$ one has
$\lambda_{t}^{n}$ is completely atomic $($or completely
non-atomic$)$ for all $n$, then so is $\lambda_{t}^{0}$ almost
everywhere.
\end{lem}
\noindent The proof is straight forward. We omit the proof. For
details check \cite{MR1000012}.

For a masa $A\subset \mathcal{M}$, fix a compact Hausdorff space $X$
such that $C(X)\subset A$ is an unital, norm separable and
\emph{w.o.t} dense $C^{*}$ subalgebra. For $\zeta\in
L^{2}(\mathcal{M})$ let $\kappa_{\zeta}:C(X)\otimes C(X)\mapsto \C$
be the linear functional defined by
\begin{align}
\nonumber \kappa_{\zeta}(a\otimes b)=\langle a\zeta b,\zeta\rangle.
\end{align}
Then $\kappa_{\zeta}$ induces an unique Radon measure $\eta_{\zeta}$
on $X\times X$ given by
\begin{align}\label{measure_from_kappa}
\kappa_{\zeta}(a\otimes b)=\int_{X\times
X}a(t)b(s)d\eta_{\zeta}(t,s)
\end{align}
and $\norm{\eta_{\zeta}}_{t.v}=\norm{\kappa_{\zeta}}$.

\indent For $\zeta_{1},\zeta_{2}\in L^{2}(\mathcal{M})$ let
$\eta_{\zeta_{1},\zeta_{2}}$ denote the possibly complex measure on
$X\times X$ obtained from the vector functional
\begin{align}\label{b011_included}
\langle a\zeta_{1} b,\zeta_{2}\rangle=\int_{X\times
X}a(t)b(s)d\eta_{\zeta_{1},\zeta_{2}}(t,s), \text{ }a,b\in C(X).
\end{align}
We will write $\eta_{\zeta,\zeta}=\eta_{\zeta}$. Note that
$\eta_{\zeta}$ is a positive measure for all $\zeta\in
L^{2}(\mathcal{M})$. It is easy to see that the following
polarization type identity holds:
\begin{align}\label{polarize}
4\eta_{\zeta_{1},\zeta_{2}}=\left(\eta_{\zeta_{1}+\zeta_{2}}-\eta_{\zeta_{1}-\zeta_{2}}\right)+i\left(\eta_{\zeta_{1}+i\zeta_{2}}-\eta_{\zeta_{1}-i\zeta_{2}}
\right).
\end{align}
Note that the decomposition of $\eta_{\zeta_{1},\zeta_{2}}$ in Eq.
\eqref{polarize} need not be its Hahn decomposition in general, but
\begin{align}\label{polarize1}
\nonumber
\abs{\eta_{\zeta_{1},\zeta_{2}}}&\leq\left(\eta_{\zeta_{1}+\zeta_{2}}+\eta_{\zeta_{1}-\zeta_{2}}\right)+\left(\eta_{\zeta_{1}+i\zeta_{2}}+\eta_{\zeta_{1}-i\zeta_{2}}\right)
=4(\eta_{\zeta_{1}}+\eta_{\zeta_{2}}).
\end{align}
So
\begin{align}
\abs{\eta_{\zeta_{1},\zeta_{2}}}\leq
\eta_{\zeta_{1}}+\eta_{\zeta_{2}}.
\end{align}

\begin{lem}\label{uniform_convergence_of_vector_functionals}
If $\zeta_{n},\zeta\in L^{2}(\mathcal{M})$ be
such that, 
$\zeta_{n}\rightarrow \zeta$ in $\norm{\cdot}_{2}$ then 
\begin{align}
\nonumber \eta_{\zeta_{n}}\rightarrow \eta_{\zeta} \text{ in }
\norm{\cdot}_{t.v}.
\end{align}
\end{lem}
\begin{proof}
Obvious.
\end{proof}

\begin{prop}\label{finite_atomic_measures}
Let $A\subset \mathcal{M}$ be a masa. Let $X$ be a compact Hausdorff
space such that $C(X)\subset A$ is unital, norm separable and w.o.t
dense in $A$ and let $\nu$ be the tracial measure. Let $0\neq
\zeta\in L^{2}(N(A)^{\prime\prime})$. Then
$\eta_{\zeta_{t}},\eta_{\zeta_{s}}$ is completely atomic $\nu$
almost all $t,s$ where $\eta_{\zeta}$ is the measure defined in Eq.
\eqref{measure_from_kappa} and $\eta_{\zeta_{t}},\eta_{\zeta_{s}}$
are $(\pi_{1},\nu)$ and $(\pi_{2},\nu)$ disintegrations of
$\eta_{\zeta}$ respectively.
\end{prop}
\begin{proof}
We only prove for the $(\pi_{1},\nu)$ disintegration. If $\zeta=u$
where $u\in N(A)$ then the result is obvious as the measure
$\eta_{u}$ will be concentrated on the automorphism graph. The span
of $N(A)$ being \emph{s.o.t} dense in $N(A)^{\prime\prime}$ it
suffices by Lemma \ref{uniform_convergence_of_vector_functionals}
and \ref{close_measure_lemma} to prove the statement when
$\zeta=\sum_{i=1}^{n}c_{i}u_{i}$ where $u_{i}\in N(A)$ and $c_{i}\in
\C$ for $1\leq i\leq n$. Now for $a,b\in A$
\begin{align}
\nonumber\langle
a(\sum_{i=1}^{n}c_{i}u_{i})b,(\sum_{i=1}^{n}c_{i}u_{i})\rangle&=\sum_{i=1}^{n}\abs{c_{i}}^{2}\langle
au_{i}b,u_{i}\rangle+\sum_{i\neq j=1}^{n}c_{i}\bar{c_{j}}\langle
au_{i}b,u_{j}\rangle.
\end{align}
The measures given by $a\otimes b\mapsto \abs{c_{i}}^{2}\langle
au_{i}b,u_{i}\rangle$, $a,b\in C(X)$ are concentrated on the
automorphism graphs implemented by $u_{i}$ and hence definitely
disintegrates as atomic measures and so does their sum from Lemma
\ref{fibre_addition_formulae}. The measures given by $a\otimes
b\mapsto c_{i}\bar{c_{j}}\langle au_{i}b,u_{j}\rangle$, $a,b\in
C(X)$ for $i\neq j$ are possibly complex measures. However Eq.
\eqref{polarize1} forces that these measures are also concentrated
on the union of the automorphism graphs implemented by $u_{i}$ and
$u_{j}$. Thus $\eta_{\sum_{i=1}^{n}c_{i}u_{i}}$ is concentrated on
the union of the automorphism graphs implemented by $u_{i}$, $1\leq
i\leq n$. Hence the result follows.
\end{proof}

\section{Fundamental Set and Generalized Dye's Theorem}

This section is intended to characterize some operators in the
normalizing algebra of a masa. Throughout this section $\mathcal{N}$
will denote a finite von Neumann algebra gifted with a faithful,
normal, normalized trace $\tau$. $B\subset \mathcal{N}$ will denote
a von Neumann subalgebra of $\mathcal{N}$.\\  
\indent As usual $\mathcal{N}$ will be assumed to be acting on
$L^{2}(\mathcal{N},\tau)$ by left multipliers.
$L^{2}(\mathcal{N},\tau)$ is a $B$-$B$ Hilbert $w^{*}$-bimodule for
any von Neumann subalgebra $B\subset \mathcal{N}$. We know if
$\mathbb{E}_{B}$ denotes the unique trace preserving conditional
expectation onto $B$, then $\mathbb{E}_{B}$ is given by the Jones
projection $e_{B}$ associated to $B$ via the formula
$\mathbb{E}_{B}(x)\hat{1}=e_{B}(x\hat{1})$. For $b_{1},b_{2}\in B$
and $\zeta\in L^{2}(\mathcal{N},\tau)$ one has
\begin{align}
e_{B}(b_{1}\zeta b_{2})=b_{1}e_{B}(\zeta)b_{2}.
\end{align}

We will interchangeably use the symbols $\mathbb{E}_{B}$ and
$e_{B}$.
\begin{defn}
For a subalgebra $B\subset \mathcal{N}$ define the \emph{fundamental
set} of $B$ to be
$$N^{f}(B)=\{x\in \mathcal{N}: Bx=xB\}.$$
\end{defn}
\indent Note that $x\in N^{f}(B)$ implies $x^{*}\in N^{f}(B)$.
\begin{defn}
For a subalgebra $B\subset \mathcal{N}$ define the
\emph{weak-fundamental set} of $B$ to be
$$N^{f}_{2}(B)=\{\zeta\in L^{2}(\mathcal{N},\tau): B\zeta=\zeta B\}.$$
\end{defn}
\indent Note that $\zeta\in N^{f}_{2}(B)$ implies $\zeta^{*}\in
N^{f}_{2}(B)$ and $N^{f}(B)\subset N^{f}_{2}(B)$. When $B$ is a
masa,
$\zeta\in N^{f}_{2}(B)$ implies $a\zeta,\zeta a\in N^{f}_{2}(B)$ for all $a\in B$. \\
\indent To understand the normaliser of a masa the set
$N^{f}_{2}(B)$ will naturally arise into the scene. However working
with vectors in $L^{2}(\mathcal{N},\tau)$ is always a technical
issue. Polar decomposition of vectors and the theory of $L^{1}$
spaces are the tools we need, for which we will give a short
exposition. For details check Appendix B of \cite{MR999996} and
\cite{MR1000012}. To keep it short we will omit most proofs. It is
here, where one usually encounters unbounded operators. For results
proved
in this section we have borrowed ideas from Roger Smith.\\
\indent The positive cone $L^{2}(\mathcal{N},\tau)^{+}$ in
$L^{2}(\mathcal{N},\tau)$ is defined to be
$\overline{{\mathcal{N}}^{+}}^{\norm{\cdot}_{2}}$ i.e. the closure
of the positive elements of $\mathcal{N}$ in
$L^{2}(\mathcal{N},\tau)$. It can be shown that
$L^{2}(\mathcal{N},\tau)$ is the algebraic span of
$L^{2}(\mathcal{N},\tau)^{+}$. For $x\in \mathcal{N}$ the equation
$\norm{x}_{1}=\tau(\abs{x})$ defines a norm on $\mathcal{N}$. The
completion of $\mathcal{N}$ with respect to $\norm{\cdot}_{1}$ is
denoted by $L^{1}(\mathcal{N},\tau)$. It can be shown that
\begin{align}\label{l_1_norm}
\norm{x}_{1}=\sup\{\abs{\tau(xy)}: y\in \mathcal{N}, \norm{y}\leq
1\}.
\end{align}
So $\abs{\tau(x)}\leq \norm{x}_{1}$. Thus by density of
$\mathcal{N}$ in $L^{1}(\mathcal{N},\tau)$, $\tau$ extends to a
bounded linear functional on $L^{1}(\mathcal{N},\tau)$ which will
also be denoted by $\tau$. One can analogously define the positive
cone of $L^{1}(\mathcal{N},\tau)$ which we denote by
$L^{1}(\mathcal{N},\tau)^{+}$. Clearly
$\norm{x}_{1}=\norm{x^{*}}_{1}$. Consequently, the Tomita operator
$J$ extends to a surjective anti-linear isometry to
$L^{1}(\mathcal{N},\tau)$ which will also be denoted by $J$.
Moreover $J^{2}=1$. We will interchangeably use
the notations $J\zeta$ and $\zeta^{*}$ for $\zeta\in L^{1}(\mathcal{N},\tau)$. \\
\indent Both the spaces $L^{1}(\mathcal{N},\tau)$ and
$L^{2}(\mathcal{N},\tau)$ are unitary $\mathcal{N}$-$\mathcal{N}$
bimodules. The space $L^{1}(\mathcal{N},\tau)$ can be identified
with the predual of $\mathcal{N}$ and $L^{2}(\mathcal{N},\tau)$ is
dense in $L^{1}(\mathcal{N},\tau)$. One also has
$\tau(x\zeta)=\tau(\zeta x)$ for $x\in \mathcal{N}$ and $\zeta\in
L^{1}(\mathcal{N},\tau)$. Note that $\mathbb{E}_{B}$ is a
contraction from $\mathcal{N}$ onto $B$. It can be shown that for
$x\in \mathcal{N}$,
\begin{align}\label{conditional_expecttaion_for_subalgebra}
\norm{\mathbb{E}_{B}(x)}_{1}\leq \norm{x}_{1}.
\end{align}
Thus $\mathbb{E}_{B}$ has an unique bounded extension to a
contraction from $L^{1}(\mathcal{N},\tau)$ onto $L^{1}(B,\tau)$,
which will as well be denoted by $\mathbb{E}_{B}$. This extension
preserves the extension of the trace $\tau$, is $B$ modular,
positive and faithful. The bilinear map $\Psi : \mathcal{N}\times
\mathcal{N}\mapsto \mathcal{N}$ defined by $\Psi(x,y)=xy$ satisfies
\begin{align}\label{product_estimate}
\norm{\Psi(x,y)}_{1} \leq \norm{x}_{2}\norm{y}_{2}
\end{align}
by Cauchy-Schwarz inequality. Therefore $\Psi$ lifts to a jointly
continuous map from $L^{2}(\mathcal{N},\tau)\times
L^{2}(\mathcal{N},\tau)$ into $L^{1}(\mathcal{N},\tau)$. The
extension is actually a surjection. Since $\Psi$ is the product map
of operators at the level of von Neumann algebra one calls
$\Psi(\zeta_{1},\zeta_{2})$ to be $\zeta_{1}\zeta_{2}$, for
$\zeta_{1},\zeta_{2}\in L^{2}(\mathcal{N},\tau)$.
\begin{lem}$($B.5.1, \cite{MR999996}$)$\label{estimate_of_norm_1_and_2}
Let $a,b\in \mathcal{N}$ be positives. Then
\begin{align}
\norm{a^{\frac{1}{2}}-b^{\frac{1}{2}}}_{2}^{2}\leq 2 \norm{a-b}_{1}.
\end{align}
\end{lem}
Elements of $L^{1}(\mathcal{N},\tau)$ and $L^{2}(\mathcal{N},\tau)$
can be regarded as unbounded operators on $L^{2}(\mathcal{N},\tau)$.
By using the unbounded operator theory for operators affiliated to
$\mathcal{N}$, for each $\zeta\in L^{1}(\mathcal{N},\tau)^{+}$ there
exists an unique $0\leq \zeta_{0}\in L^{2}(\mathcal{N},\tau)$ such
that $\zeta_{0}^{*}\zeta_{0}=\zeta_{0}^{2}=\zeta$. In this case,
$\zeta_{0}$ is said to be the \emph{square root} of $\zeta$ and one
writes $\zeta_{0}=\sqrt{\zeta}=\zeta^{\frac{1}{2}}$. For $\zeta\in
L^{2}(\mathcal{N},\tau)$ one has $\zeta^{*}\zeta\in
L^{1}(\mathcal{N},\tau)$. From Eq. \ref{product_estimate} and Lemma
\ref{estimate_of_norm_1_and_2} it follows that $\zeta^{*}\zeta\in
L^{1}(\mathcal{N},\tau)^{+}$. In particular,
$\sqrt{\zeta^{*}\zeta}\in L^{2}(\mathcal{N},\tau)$ for any $\zeta\in
L^{2}(\mathcal{N},\tau)$ and the square root of any positive in
$L^{1}(\mathcal{N},\tau)$ is an unique element of
$L^{2}(\mathcal{N},\tau)$. One also writes
$\abs{\zeta}=\sqrt{\zeta^{*}\zeta}$ for $\zeta\in
L^{2}(\mathcal{N},\tau)$. If $\zeta\in L^{1}(\mathcal{N},\tau)$ be
self adjoint i.e. $\zeta=\zeta^{*}$ then $\zeta=\zeta_{+}-\zeta_{-}$
where $\zeta_{\pm}\in L^{1}(\mathcal{N},\tau)^{+}$ and this
decomposition is unique by requiring that
${\zeta_{+}}^{\frac{1}{2}}{\zeta_{-}}^{\frac{1}{2}}=0$.\\
\indent Let $\zeta\in L^{2}(\mathcal{N},\tau)$. Consider the
projections $p,q$ in $\textbf{B}(L^{2}(\mathcal{N},\tau))$ whose
ranges are $\overline{J\mathcal{N}J\sqrt{\zeta^{*}\zeta}}$,
$\overline{J\mathcal{N}J\zeta}$ respectively. Since the ranges of
$p,q$ are invariant subspaces of
$J\mathcal{N}J=\mathcal{N}^{\prime}$ so $p,q$ lies in $\mathcal{N}$.
Using unbounded operators one obtains \emph{polar decomposition} of
vectors $($Eq. \eqref{polar_decopm_starts_here}$)$ which we
formalize below.
\begin{thm}\label{theorem_for_weak_polar_decompose}
There is an unique partial isometry $v\in \mathcal{N}$ with initial
projection $p$ and final projection $q$ which satisfy the following
condition:
\begin{align}\label{limit_partial_isometry_defined}
vJx^{*}J\sqrt{\zeta^{*}\zeta}=Jx^{*}J\zeta, \text{ }x\in
\mathcal{N}.
\end{align}
In particular,
\begin{align}\label{polar_decopm_starts_here}
v\sqrt{\zeta^{*}\zeta}=\zeta.
\end{align}
$(i)$ Let $B\subset
\mathcal{N}$ be a masa, then $\zeta\in L^{2}(B,\tau)$ imply $p,q\in B$.\\
$(ii)$ For $\zeta\in L^{2}(\mathcal{N},\tau)$ if $\zeta^{*}\zeta\in
\mathcal{N}$ then $\zeta\in \mathcal{N}$.
\end{thm}

For $\zeta\in L^{2}(\mathcal{N},\tau)$ we define the \emph{left and
right kernel} of $\zeta$ to be respectively $Ker_{l}(\zeta)=\{x\in
\mathcal{N}:\zeta x=0\}$ and $Ker_{r}(\zeta)=\{x\in
\mathcal{N}:x\zeta =0\}$. Then $Ker_{l}(\cdot),Ker_{r}(\cdot)$
are subspaces of $\mathcal{N}$. $Ker_{l}(\cdot)$, $Ker_{r}(\cdot)$ are $w.o.t$ and $s.o.t$ closed.\\
\indent If $\zeta\in L^{1}(\mathcal{N},\tau)$ then \emph{the}
\emph{left and the right kernels} of $\zeta$ can be defined
analogously. We will denote the \emph{kernels} of the $L^{1}$
vectors by $Ker_{l}(\cdot),Ker_{r}(\cdot)$ as well. This is slight
abuse of notation. In this case, they
are norm closed subspaces of $\mathcal{N}$.\\
\indent For $\zeta\in L^{2}(\mathcal{N},\tau)$ we have
\begin{align}\label{kernel_wot_closed}
Ker_{l}(\zeta)=Ker_{l}(\sqrt{\zeta^{*}\zeta})=Ker_{l}(\zeta^{*}\zeta).
\end{align}
However the righthand side is defined in $L^{1}$ sense. Therefore
for $\zeta\in L^{2}(\mathcal{N},\tau)$, $Ker_{l}(\zeta^{*}\zeta)$
$($respectively $Ker_{r}(\zeta\zeta^{*})$$)$ are in fact $w.o.t$
closed. Similar statements hold for $Ker_{r}(\cdot)$ as well.\\
\indent For $\zeta\in L^{2}(\mathcal{N},\tau)$ we define the
\emph{left and right ranges} of $\zeta$ to be respectively
$Ran_{l}(\zeta)=\{\zeta x: x\in \mathcal{N}\}$ and
$Ran_{r}(\zeta)=\{x\zeta : x\in \mathcal{N}\}$.\\
Note that for $\zeta\in L^{2}(\mathcal{N},\tau)$,
\begin{align}\label{kernel_and_range_relation}
 \{x\in \mathcal{N}:\zeta x=0\}&=\{x\in \mathcal{N}:\langle \zeta x,y\rangle=0 \text{ for all }y\in \mathcal{N}\}\\
\nonumber &=\{x\in \mathcal{N}:\langle x,\zeta^{*} y\rangle=0 \text{
for all }y\in \mathcal{N}\}
\end{align}
implies $Ker_{l}(\zeta)=Ran_{l}(\zeta^{*})^{\perp}$.

\begin{prop}\label{identify_initial_and_final_ranges}
Let $\zeta\in L^{2}(\mathcal{N},\tau)$ and let
$\zeta=v\sqrt{\zeta^{*}\zeta}$ be its polar decomposition. Then
$v^{*}v$ is the projection from $L^{2}(\mathcal{N},\tau)$ onto
$Ker_{l}(\zeta)^{\perp}$ and $vv^{*}$ is the projection onto
$\overline{Ran_{l}(\zeta)}$.
\end{prop}

\begin{prop}\label{identify_the_range_proj_as_limit_anywhere}
Let $\zeta\in L^{2}(\mathcal{N},\tau)$ and let $\zeta=v\abs{\zeta}$
be its polar decomposition. Then
$\abs{\zeta}^{\frac{1}{2^{k}}}\rightarrow v^{*}v$ as $k\rightarrow
\infty$ in $\norm{\cdot}_{2}$.
\end{prop}
The proof of Prop. \ref{identify_the_range_proj_as_limit_anywhere}
is a direct application of monotone convergence theorem.

\begin{lem}\label{l^{1}A_vectors}
Let $A\subset \mathcal{N}$ be a masa. Let $\zeta\in
L^{1}(\mathcal{N},\tau)$ be a nonzero vector such that $a\zeta=\zeta
a$ for all $a\in A$. Then $\zeta\in L^{1}(A,\tau)$.
\end{lem}
\begin{proof}
First assume $\zeta\geq 0$. Then use uniqueness of square roots of
$L^{1}$ vectors. In, the general case write $\zeta$ as a linear
combination of four positives. We omit the details.
\end{proof}

\begin{prop}\label{at_the_level_or_one_sp_normaliser}
Let $A\subset \mathcal{N}$ be a masa. Let $0\neq\zeta\in
L^{1}(\mathcal{N},\tau)^{+}$ be such that $A\zeta=\zeta A$. Then
$\zeta\in {L^{1}(A,\tau)}^{+}$.
\end{prop}
\begin{proof}
Let $\mathcal{I}=\{a\in A: a\zeta=0\}$. Then $\mathcal{I}$ is a
weakly closed ideal $($see Eq. \eqref{kernel_wot_closed} and related
discussion$)$ in $A$ and so has the form $A(1-p)$ for some
projection $p\in A$. Then $p\zeta=\zeta$, so $\zeta=\zeta p$ by
operating with extended
Tomita's involution operator. Thus $Ap\zeta=A\zeta p=\zeta Ap$.\\
For $a_{1},a_{2}\in A$ if $\zeta a_{1}p=\zeta a_{2}p$ then
$\zeta(a_{1}-a_{2})p=0$, so $p(a_{1}^{*}-a_{2}^{*})\zeta=0$. Hence
$p(a_{1}^{*}-a_{2}^{*})\in \mathcal{I}$, but $1-p$ is the identity
for $\mathcal{I}$. So $p(a_{1}^{*}-a_{2}^{*})=0$ and
hence $a_{1}p=a_{2}p$.\\
This means there is a well defined map $\psi :Ap\mapsto Ap$ such
that
$$ap\zeta=\zeta\psi(ap) \text{ for } a\in A.$$ Taking conditional
expectation $($see Eq.
\eqref{conditional_expecttaion_for_subalgebra} and related
discussion$)$ one gets $(ap-\psi(ap))\mathbb{E}_{A}(\zeta)=0$ $($the
left and the right action by elements of $A$ coincides on
$L^{1}(A,\tau)$$)$. Suppose there is an operator $a\in A$ such that
$ap-\psi(ap)\neq 0$. Write $ap-\psi(ap)=bp$ for $b\in A$. Then
$pb^{*}bp\mathbb{E}_{A}(\zeta)=0$, so
$\mathbb{E}_{A}(pb^{*}bp\zeta)=0$. 
Let
$\zeta=\underset{n}\lim\text{ } x_{n}$ in $\norm{\cdot}_{1}$ where
$x_{n}\in \mathcal{N}^{+}$. Therefore
$$\underset{n}\lim\text{
}
\tau(x^{\frac{1}{2}}_{n}(bp)^{*}bpx_{n}^{\frac{1}{2}})=\underset{n}\lim\text{
} \tau(pb^{*}bpx_{n})=\underset{n}\lim\text{ }
\tau(\mathbb{E}_{A}(pb^{*}bpx_{n}))=0.$$ The last statement follows
from Eq. \eqref{l_1_norm} and Eq.
\eqref{conditional_expecttaion_for_subalgebra}. 
So $\underset{n}\lim\text{ }bpx^{\frac{1}{2}}_{n}=0$ in
$\norm{\cdot}_{2}$ and hence $bp\zeta=\underset{n}\lim\text{
}bpx_{n}=0$, in $\norm{\cdot}_{1}$ by Lemma
\ref{estimate_of_norm_1_and_2} and Eq. \eqref{product_estimate}.
Thus $bp\in \mathcal{I}$ so $bp=bp(1-p)=0$, a contradiction. Thus
$\psi(ap)=ap$ for all $a\in
A$.\\
Now $\zeta\in L^{1}(p\mathcal{N}p,\tau)$ and $Ap$ is a masa in
$p\mathcal{N}p$, thus $\zeta\in L^{1}(Ap,\tau)$ as
$ap\zeta=\zeta\psi(ap)=\zeta ap$ for all $a\in A$, from Lemma
\ref{l^{1}A_vectors}.
\end{proof}

\begin{thm}$($Generalized Dye's theorem-$L^{2}$ form$)$\label{identify_operators_in_q_Hilbert_normalisers}
Let $A\subset \mathcal{N}$ be a masa. Then $\zeta\in N^{f}_{2}(A)$
if and only if $\zeta=v\xi$ for some $\xi\in L^{2}(A,\tau)$ and
$v\in \mathcal{GN}(A)$. In particular, $\overline{span
N^{f}(A)}^{\norm{\cdot}_{2}}=L^{2}(N(A)^{\prime\prime},\tau)$.
\end{thm}
\begin{proof}
\textbf{Case 1}: Assume $\zeta\in N^{f}_{2}(A)$ and $\zeta\geq 0$
i.e. $\zeta\in {\overline{\mathcal{N}^{+}}}^{\norm{\cdot}_{2}}$.
Then $\zeta\in L^{1}(\mathcal{N},\tau)^{+}$ as well. From Prop.
\ref{at_the_level_or_one_sp_normaliser} we get $\zeta\in
L^{1}(A,\tau)\cap L^{2}(\mathcal{N},\tau)=L^{2}(A,\tau)$.\\
\textbf{Case 2}: Let $\zeta\in N^{f}_{2}(A)$. We may without loss of
generality assume that $\norm{\zeta}_{2}=1$. Then as $A\zeta=\zeta
A$ we also have $A\zeta^{*}=\zeta^{*}A$. So
$A\zeta^{*}\zeta=\zeta^{*}A\zeta =\zeta^{*}\zeta A$. From Prop.
\ref{at_the_level_or_one_sp_normaliser},
$$\zeta^{*}\zeta\in L^{1}(A,\tau)$$ and similarly we have
$\zeta\zeta^{*}\in L^{1}(A,\tau)$. Then
$\norm{\zeta^{*}\zeta}_{1}\leq 1$.\\
Arguing as in Prop. \ref{at_the_level_or_one_sp_normaliser}, there
are projections $p_{1},p_{2}\in A$ such that $J_{1}=\{a\in
A:a\zeta=0\}=A(1-p_{1})$ and $J_{2}=\{a\in
A:\zeta a=0\}=A(1-p_{2})$. Therefore we have $p_{1}\zeta=\zeta$ and $\zeta p_{2}=\zeta$.\\
Then there is a well defined map $($as explained before$)$
$\psi:Ap_{1}\mapsto Ap_{2}$ such that
$$ap_{1}\zeta=\zeta\psi(ap_{1})\text{ for all }a\in A.$$
Let $\zeta=v\sqrt{\zeta^{*}\zeta}$ be the polar decomposition of
$\zeta$ from Thm. \ref{theorem_for_weak_polar_decompose}. Then $v$
is a partial isometry in $\mathcal{N}$ and the initial space of $v$
is\\ ${\{\sqrt{\zeta^{*}\zeta}x: x\in
\mathcal{N}\}}^{-\norm{\cdot}_{2}}$ and the final space is $\{\zeta
x:x\in \mathcal{N}\}^{-\norm{\cdot}_{2}}$. Moreover
the projections $v^{*}v$ and $vv^{*}$ are in $A$.\\
Indeed, by Prop. \ref{identify_initial_and_final_ranges}, $v^{*}v$
is the projection onto $Ker_{l}(\zeta)^{\perp}$ and $vv^{*}$ onto
$\overline{Ran_{l}(\zeta)}$. By Prop.
\ref{identify_the_range_proj_as_limit_anywhere}, $v^{*}v\in A$.
Replacing $\zeta$ by $\zeta^{*}$ and using
$Ker_{l}(\zeta)^{\perp}=\overline{Ran_{l}(\zeta^{*})}$ $($see Eq.
\eqref{kernel_and_range_relation}$)$, a similar argument will yield
$vv^{*}\in A$. Clearly $v^{*}v=p_{2}$ and $vv^{*}=p_{1}$. Then
$$ap_{1}v\sqrt{\zeta^{*}\zeta}=v\sqrt{\zeta^{*}\zeta}\psi(ap_{1}).$$
Now $$J_{0}=\{b\in A:ap_{1}vb=vb\psi(ap_{1})\text{ for all }a\in
A\}$$ is a weakly closed ideal in $A$ and its closure in
$\norm{\cdot}_{2}$ is precisely the set
$$J_{0}^{-\norm{\cdot}_{2}}=\{\xi\in
L^{2}(A,\tau):ap_{1}v\xi=v\xi\psi(ap_{1})\text{ for all }a\in A\}$$
which
contains $\sqrt{\zeta^{*}\zeta}$.\\
Since the left and right action of $A$ on $L^{2}(A,\tau)$ agree, so
$\xi_{0}\in J_{0}^{-\norm{\cdot}_{2}}$ and $a\in A$ implies that
$\xi_{0}
a,a\xi_{0}\in J_{0}^{-\norm{\cdot}_{2}}$.\\
Since the $w.o.t$ closed ideal $J_{0}$ in $A$ is just a cutdown of
$A$ by a projection from $A$ any positive $\zeta_{0}\in
J_{0}^{-\norm{\cdot}_{2}}$ is a limit in $\norm{\cdot}_{2}$ of an
increasing sequence of positive operators from $J_{0}$. Now it
follows that $\abs{\zeta}^{\frac{1}{2^{k}}}\in
J_{0}^{-\norm{\cdot}_{2}}$ for all $k\in \mathbb{N}$. Therefore by
Prop. \ref{identify_the_range_proj_as_limit_anywhere} it follows
that $v^{*}v=p_{2}\in J_{0}^{-\norm{\cdot}_{2}}$ and hence $p_{2}\in
J_{0}\subseteq A$. Similarly arguing with $\zeta\zeta^{*}$ one shows
$p_{1}\in A$. Therefore
$$ap_{1}vp_{2}=vp_{2}\psi(ap_{1})\text{ for all }a\in A.$$ Then
$$v^*av=(vp_{2})^{*}avp_{2}=v^{*}ap_{1}v=v^{*}vp_{2}\psi(ap_{1})=\psi(ap_{1}).$$
Therefore $v^{*}$ and hence $v$ are groupoid normalisers. So
$$\zeta=v\xi.$$ for $v\in \mathcal{GN}(A)$ and $\xi=\abs{\zeta}\in
L^{2}(A,\tau)^{+}$.
\end{proof}


\section{Characterization by Baire Category Methods}

The study of Cartan masas in $\rm{II}_{1}$ factors has received
special attention by many experts. Our approach of studying
\emph{measure-multiplicity-invariant} was also considered implicitly
by Popa and Shlyakhtenko in \cite{MR1951447}. 
In this
section we will use an alternative approach to characterize masas by
their \emph{left-right-measure}. As it turns out, many known
theorems related to structure and normalisers of masas that were
solved using different
techniques can be solved by a single technique. \\
\indent Let $A=L^{\infty}(X,\nu_{X}),B=L^{\infty}(Y,\nu_{Y})$ be two
diffuse commutative von Neumann algebras, where $\nu_{X},\nu_{Y}$
are probability measures. Let $C(A,B)$ denote the set of all
$A,B$-bimodules. This set $C(A,B)$ contains three distinguished
subsets.\\
\indent We will use the variable $s$ to denote the first variable
and $t$ to denote the second variable. Following \cite{MR1951447} we
define:
\begin{defn}\label{discrete_bimodule}
A \emph{discrete} $($respectively, \emph{diffuse}, \emph{mixed}$)$
$A,B$-bimodule is a Hilbert space $\h$ so that $\h\cong
\underset{i\in I}\oplus L^{2}(X\times Y,\mu_{i})$ where for all $i$,
$\mu_{i}$ disintegrates as $\mu_{i}(s,t)=\mu_{t}^{(i)}(s)\nu_{Y}(t)$
with $\mu_{t}^{(i)}$ atomic $($respectively non-atomic, a
combination of both nonzero atomic part and nonzero non-atomic
part$)$ for $\nu_{Y}$ almost all $t$.
\end{defn}
\indent It is to be noted that in view of Lemma
\ref{equivalence_of_measure_imply_equivalence_of_fibre_almost_everywhere},
the definition above only cares about the equivalence class of the
measures $\mu_{i}$ and not a particular member of the class. The
definition forces $\mu_{i}$ to be a non-atomic measure, and the
existence of such a disintegration actually forces the push forward
of $\mu_{i}$'s on the space $Y$ to be dominated by $\nu_{Y}$. We
will restrict ourselves to the case $I$ is countable. Let
$C_{d}(A,B),C_{n.a}(A,B),C_{m}(A,B)$ denote the set of
all discrete, diffuse, mixed $A,B$-bimodules respectively.\\
\indent Denote by $C_{d}(A)\subset C_{d}(A,A)\subset C(A,A)$ the set
of those bimodules $\h\in C_{d}(A,A)$ for which $\bar{\h}\in
C_{d}(A,A)$. Here $\bar{\h}$ is the opposite Hilbert space of $\h$
with left and right actions interchanged. Bimodules in $C_{d}(A)$
are precisely those for which the associated measures $\mu_{i}$'s in
Defn. \ref{discrete_bimodule} also have a completely atomic
$\nu_{X}$ disintegration. Similarly define $C_{n.a}(A),C_{m}(A)$.
Note that the spaces $C_{d}(A),C_{n.a}(A),C_{m}(A)$ are all closed
with respect to taking sub bimodules.\\ 
\indent When $A,B$ are masas in a $\rm{II}_{1}$ factor $\mathcal{M}$
the standard Hilbert space $L^{2}(\mathcal{M})$ is naturally a
$w^{*}$-continuous $A,B$ bimodule, meaning it carries a pair of
mutually commuting normal representations of $A$ and $B$.\\
\indent Note that when we deal with the \emph{left-right-measure} of
a masa, knowing the disintegration along the second variable enables
us to know the disintegration along the first variable as well, by
pushing forward the former with the flip map $($see Lemma
\ref{atom_shift_lemma_under_flip_map}$)$.\\
\indent Before we proceed to the characterization of masas we will
have to make few definitions and statements that are very valuable
tools yet not appear in standard measure theory courses. For details
see \cite{MR1725642}, \cite{MR0393403}.

\begin{defn}
Let $X$ be a Polish space. A subset $B$ of $X$ is said to have
\emph{Baire property} if there is an open set $\mathcal{O}\subset X$
and a comeager set $A\subset X$ such that $A\cap \mathcal{O}=A\cap
B$.
\end{defn}
The collection of sets with \emph{Baire property} forms a
$\sigma$-algebra which 
includes the \emph{Borel} $\sigma$-\emph{algebra}. 

\begin{defn}
Let $X$ and $Y$ be Polish spaces. A function $f:X\mapsto Y$ is said
to be \emph{Baire measurable} if the inverse image of any \emph{open
set} has \emph{Baire property}. The function $f$ is said to be
\emph{universally Baire measurable} if given any Borel function $g$
into $X$ the function $f\circ g$ is Baire measurable.
\end{defn}
\indent Note that in particular every Borel function is \emph{Baire
measurable}. 

\begin{defn}
A subset $E$ 
of a Polish space is said to be \emph{universally measurable} if it
is measurable with respect to any \emph{complete Borel probability
measure}.
\end{defn}
\begin{defn}\label{analytic_sets}
A subset $E$ of a Polish space $X$ is said to be
${\underset{\sim}\Sigma}{ }_{1}^{1}$ or \emph{analytic}, if there is
a Polish space $Y$, a Borel subset $B$ of $Y$ and a Borel function
$f:Y\mapsto X$ such that $f(B)=E$. In other words,
${\underset{\sim}\Sigma}{ }_{1}^{1}$ sets are Borel images of Borel
sets.
\end{defn}
\begin{rem}
The above definition of analytic sets is as per \cite{MR1725642}.
However in, \cite{MR1468230} 
continuous images
rather than Borel images are used. The two definitions are in fact equivalent.
\end{rem}
\indent A very nontrivial theorem of Lusin says the following.
\begin{thm}$($Lusin$)$\label{baire_lusin}
Every ${\underset{\sim}\Sigma}{ }_{1}^{1}$ set has Baire property.
Every ${\underset{\sim}\Sigma}{ }_{1}^{1}$ set is universally
measurable.
\end{thm}

For a function $f:Y\mapsto X$, the graph of $f$ will be denoted by
$\Gamma(f)=\{(f(y),y):y\in Y\}$. The next theorem is very crucial in
all our analysis.

\begin{thm}$($Selection Principle - Jankov, von Neumann$)$\label{Selection_Principle} Let $X,Y$
be Polish spaces and let $E\subset X\times Y$ be in
${\underset{\sim}\Sigma}{ }_{1}^{1}$. Then $E$ can be uniformized by
a function that is both Baire and universally measurable, in the
sense that for some $h:Y\mapsto X$ we have 
$$\Gamma (h_{\mid \pi_{Y}(E)})\subseteq E$$
with the property that $h^{-1}(U)$ has the Baire property and is
measurable with respect to any Borel probability measure for all
open $U\subseteq X$.
\end{thm}

\begin{rem}\label{h_is_measurable}
Let $\nu_{X}$ and $\nu_{Y}$ be any two Borel probability measures on
$X,Y$ respectively. Let $\sigma_{\nu_{X}}$ and $\sigma_{\nu_{Y}}$ be
the $\sigma$-algebras associated to the measures $\nu_{X}$,
$\nu_{Y}$ respectively. If $h$ is the function in Thm.
\ref{Selection_Principle}, then the inverse image of any Borel set
in $X$ under $h$ will lie in $\sigma_{\nu_{Y}}$, because the
collection of subsets of $X$ whose inverse images fall in
$\sigma_{\nu_{Y}}$ is a $\sigma$-algebra and contains all open sets.
If in addition, $h$ satisfies the property that $\nu_{X}(h(F))=0$ if
and only if $\nu_{Y}(F)=0$, then $h$ is
$(\sigma_{\nu_{Y}},\sigma_{\nu_{X}})$ measurable.
\end{rem}

Let $A\subset \mathcal{M}$ be a masa. Without loss of generality we
assume that $A=L^{\infty}([0,1],\lambda)$ where $\lambda$ is the
Lebesgue measure on $[0,1]$. Let $[\eta_{[0,1]\times [0,1]}]$ denote
the \emph{left-right-measure} of $A$. We are including the diagonal.
Fix any member $\eta_{[0,1]\times [0,1]}$ from the equivalence
class. 
Since our base space is now
fixed we will rename $\eta_{[0,1]\times [0,1]}$ by $\eta$ to reduce the notation. We assume that $\eta$ is a finite measure.\\
\indent Consider the set $S_{a}=([0,1]\times [0,1])_{a}$ as defined
in Prop. \ref{atomic_part_measurable_set} with respect to the
disintegration along the $y$-axis i.e. the $t$ variable. Then by
Prop. \ref{atomic_part_measurable_set}, $S_{a}$ is a
$[\eta]$-measurable set, i.e. measurable with respect to the
completion $\sigma$-algebra associated to $\eta$. Define measures
$$\eta_{a}=\eta_{\mid (S_{a}\setminus \Delta([0,1]))}\text{ and
}\eta_{n.a}=\eta_{\mid (S_{a}^{c}\setminus \Delta([0,1]))}.$$ Then\\
$(i)$ $\eta_{\mid \Delta([0,1])^{c}}=\eta_{a}+\eta_{n.a}, \text{ }\eta_{a}\perp \eta_{n.a}.$ \\
$(ii)$ Both $\eta_{a},\eta_{n.a}$ have disintegrations along the
$x,y$ axes with respect to $\lambda$.\\
\indent Note that the disintegration of the measure $\eta_{a}$ along
the $x$ and $y$-axes must have at most countably many atoms almost
all fibres $($see Lemma \ref{atom_shift_lemma_under_flip_map}$)$,
otherwise $\eta$ is an infinite measure. Since changing the measure
$\eta_{a}$ or $\eta$ on a set of measure $0$ does not change the
measure class of $\eta_{a}$ or $\eta$, we can as well assume without
loss of generality that, the disintegration of the measure
$\eta_{a}$ along $y$-axis $($second variable$)$ has at most
countable number of atoms for all fibres. With this as set up we
formalize the main theorem of this manuscript. Thm.
\ref{classification_theorem} will be proved latter in this section.

\begin{thm}$($Classification of Types$)$\label{classification_theorem}
A masa $A\subset \mathcal{M}$ is
\begin{align}
\nonumber &(i) \text{Cartan if and only if } \eta_{n.a}=0 \text{ equivalently }L^{2}(A)^{\perp}\in C_{d}(A),\\
\nonumber &(ii) \text{singular if and only if } \eta_{a}=0 \text{ equivalently }L^{2}(A)^{\perp}\in C_{n.a}(A),\\
\nonumber &(iii) A \varsubsetneq N(A)^{\prime\prime}\varsubsetneq
\mathcal{M} \text{ if and only if }\eta_{a}\neq 0,\eta_{n.a}\neq 0
\text{ equivalently }\\
\nonumber &\indent\indent L^{2}(A)^{\perp}\in C_{m}(A).\\
\nonumber &(iv) A \text{ is semiregular if and only if the closed
support of }\eta_{a}\\
\nonumber &\indent \indent\text{is }[0,1]\times [0,1].
\end{align}
\end{thm}
\begin{rem}
First of all, in view of Lemma \ref{fibre_addition_formulae} and
\ref{equivalence_of_measure_imply_equivalence_of_fibre_almost_everywhere},
the characterization does not depend on any particular member
of the \emph{left-right-measure}.\\
Secondly, $L^{2}(A)$ is always included in $C_{d}(A)$, the
disintegration having one atom at each point of the diagonal. That
is the reason one excludes $L^{2}(A)$ from statements in Thm.
\ref{classification_theorem}.\\
Finally, from our discussion on direct integrals in Sec. $2$, it
follows that $L^{2}(A)^{\perp}$ is the direct integral over
$[0,1]\times [0,1]$ with respect to the measure $\eta_{\mid (\Delta
[0,1])^{c}}$, the measurable field of Hilbert spaces depending on
$m_{[0,1]}$ or the \emph{Puk\'{a}nszky invariant}. So the equivalent
statements regarding the type of bimodules and measure in Thm.
\ref{classification_theorem} are obvious statements.
\end{rem}

The next technical lemma is the key to characterization of masas.
There are several measures involved in its statement and proof.
Since there is danger of confusion with measurability of objects
involved we will always use phrases like ``$\mu$-measurable".

\begin{lem}\label{construction_of_Baire_functions}
Let $\eta_{a}\neq 0$. Let $Y\subseteq (\Delta[0,1])^{c}$ be a
$\eta$-measurable set of strictly positive $\eta_{a}$-measure. There
exists a $\lambda$-measurable set $E^{Y}\subseteq [0,1]$ with
$\lambda(E^{Y})>0$ and a function $h_{Y}:[0,1]\mapsto [0,1]$ such that\\
$(i)$ $h_{Y}$ is $\lambda$-measurable,\\
$(ii)$ $\Gamma(h_{Y})$ is a $\eta$-measurable set,\\
$(iii)$ $\eta(\Gamma(h_{Y}))>0$ and $(h_{Y}(t),t) \in Y\cap S_{a}$ for $t\in E^{Y}$,\\
$(iv)$ for $E\subset [0,1]$, $\lambda(E)=0$ if and only if
$\lambda(h_{Y}(E))=0$.
\end{lem}
\begin{proof}
We have
$$\eta(S_{a}\cap Y)=\eta_{a}(S_{a}\cap Y)=\eta_{a}(Y)>0.$$
Consider the disintegration of $\eta_{\mid Y}$ along the $y$-axis.
There is a set $F^{Y}\subseteq [0,1]$ such that $\lambda(F^{Y})>0$
and for each $t\in F^{Y}$ the measure $({\eta_{\mid Y}})_{t}$ has
atoms with at most countable number of atoms and for $t\not\in
F^{Y}$ the same disintegration has no atoms. This is true because
$\eta$ is a finite measure, the set $F^{Y}$ being $\pi_{y}(S_{a}\cap
Y)$, $\pi_{y}$ denoting the projection on to the $y$-axis. The set
$S_{a}\cap Y$ is $\eta$-measurable, so $S_{a}\cap Y=B\cup N$ where
$B$ is a Borel set in $[0,1]\times [0,1]$ and $N$ is a $\eta$-null
set. The set $B$ is a continuous image of a Polish space by Thm.
14.3.5 of \cite{MR1468230} and so is $\pi_{y}(B)$. By Defn
\ref{definition_of_disintegration}, $\lambda(\pi_{y}(N))=0$. So
$F^{Y}$ is $\lambda$-measurable set by Thm. 
\ref{baire_lusin}. Throwing off another $\lambda$-null set from
$F^{Y}$ if necessary we can as well assume without loss of
generality that $F^{Y}$ is a
Borel set.\\
Let 
$F_{a}^{Y}=\left((Y\cap S_{a})\cap ([0,1]\times F^{Y})\right)$ which
is $\eta$-measurable. Write $F_{a}^{Y}=E_{a}^{Y}\cup N_{1}$ where
$N_{1}$ is a $\eta$-null set and $E_{a}^{Y}$ is a Borel set.
Then by Thm. $14.3.5$ of \cite{MR1468230}, 
$E_{a}^{Y}$ is in ${\underset{\sim}\Sigma}{ }_{1}^{1}$, in fact it
is the continuous
image of a Polish space. The hypothesis guarantees $\eta(E_{a}^{Y})>0$.\\
\noindent Let $E^{Y}=\pi_{y}(E_{a}^{Y})$. Then $E^{Y}$ is in
${\underset{\sim}\Sigma}{ }_{1}^{1}$ and hence $E^{Y}$ is
$\lambda$-measurable by Thm. \ref{baire_lusin}. 
Therefore by
Def \ref{definition_of_disintegration}, $\lambda(E^{Y})>0$. By Thm.
\ref{Selection_Principle} applied to $E_{a}^{Y}$, there exists a
function $h_{Y}: [0,1]\mapsto [0,1]$ that is both \emph{Baire} and
\emph{universally measurable} in the sense of Thm.
\ref{Selection_Principle}, such that $\Gamma({h_{Y}}_{\mid
E^{Y}})\subseteq E_{a}^{Y}$.\\
%
The inverse image under $h_{Y}$ of any Borel subset of $[0,1]$
belongs to $\sigma_{\lambda}$. Therefore given $\epsilon>0$, by
Lusin's theorem there is a closed subset $G^{Y}\subseteq E^{Y}$ such
that $\lambda(E^{Y}\setminus G^{Y})<\epsilon$ and ${h_{Y}}_{\mid
G^{Y}}$ is continuous. Then ${h_{Y}}_{\mid G^{Y}}$ is Borel
measurable. So by Cor. 2.11 of \cite{MR591683},
$\Gamma({h_{Y}}_{\mid G^{Y}})$ is Borel measurable
and hence $\eta$-measurable.\\
The disintegration along the $y$-axis of the measure
$\eta_{\mid\Gamma({h_{Y}}_{\mid G^{Y}})}$ is precisely the atom at
the point $(h_{Y}(t),t)$ for each $t\in G^{Y}$ of the measure
$\eta_{t}$. Outside $G^{Y}$ we don't care. If
$\eta(\Gamma({h_{Y}}_{\mid G^{Y}}))=0$ then by definition of
disintegration
$$0=\int_{G^{Y}}\eta_{t}(\Gamma(h_{Y}))d\lambda(t)$$ which implies
that for $\lambda$ almost all $t\in G^{Y}$ the point $(h_{Y}(t),t)$
is not an atom of $\eta_{t}$ and hence cannot be in $S_{a}$.
So $\eta(\Gamma({h_{Y}}_{\mid G^{Y}}))>0$.\\
Clearly, ${h_{Y}}_{\mid G^{Y}}$ satisfies the property that for any
$E\subset G^{Y}$, $\lambda(E)=0$ if and only if
$\lambda(h_{Y}(E))=0$. Therefore by Thm.
\ref{Tidze_extension_obeying_null}, extend ${h_{Y}}_{\mid G^{Y}}$ to
a continuous function $\tilde{h}_{Y}$ which satisfies the property
that for any $E\subset [0,1]$, $\lambda(E)=0$ if and only if
$\lambda(\tilde{h}_{Y}(E))=0$. So by Rem \ref{h_is_measurable},
$\tilde{h}_{Y}$ is $(\sigma_{\lambda},\sigma_{\lambda})$ measurable.
Rename $\tilde{h}_{Y}$ to $h_{Y}$ and $G^{Y}$ to $E^{Y}$. The rest is clear from construction.
\end{proof}

\begin{lem}\label{find_partial_isometry}
Let $\eta_{a}\neq 0$. Let $Y\subseteq (\Delta[0,1])^{c}$ be a
$\eta$-measurable set of strictly positive $\eta_{a}$-measure. Then
$\mathcal{U}(A)\varsubsetneq N(A)$, where $\mathcal{U}(A)$ denotes
the unitary group of $A$.\\
More precisely, there exists a subset $F^{Y}$ of $[0,1]$ such that
$\lambda(F^{Y})>0$, a invertible map $h_{Y}: F^{Y}\mapsto
h_{Y}(F^{Y})$ and a nonzero vector $\zeta_{Y}\in
L^{2}(N(A)^{\prime\prime})\ominus L^{2}(A)$ such that
\begin{align}
\nonumber &(i) \text{ }\zeta_{Y}=v_{Y}\rho_{Y}\text{ with }v_{Y}\in \mathcal{GN}(A),\rho_{Y}\in L^{2}(A)^{+}\\
\nonumber &(ii)\text{ }\overline{A\zeta_{Y}A}^{\norm{\cdot}_{2}}\cong\int_{\Gamma(h_{Y})}^{\oplus}\C_{s,t}d\eta(s,t),\text{ where }\C_{s,t}=\C,\\
\nonumber &(iii)\text{ }\Gamma(h_{Y})\subseteq Y\cap S_{a},\\
\nonumber &(iv)\text{ }\eta(\Gamma(h_{Y}))>0,\\
\nonumber &(v)1\in Puk(A).
\end{align}
\end{lem}
\begin{proof}
Using Lemma \ref{construction_of_Baire_functions}, choose the
function $h_{Y}$ that satisfies the conclusion of that Lemma. Note
that $h_{Y}$ satisfies the conditions of Prop.
\ref{structure_of_measurable_functions}. Apply Prop.
\ref{structure_of_measurable_functions} to the function $h_{Y}$ and
the set $E^{Y}$ to extract a set $F^{Y}\subseteq E^{Y}$ such that
$\lambda(F^{Y})>0$ and $h_{Y}$ is one to one on $F^{Y}$. So
\begin{align}
\nonumber h_{Y}: F^{Y}\mapsto h_{Y}(F^{Y}) \text{ is invertible.}
\end{align}
Note that 
as $\lambda(F^{Y})>0$ so $\eta(\Gamma({h_{Y}}_{\mid F^{Y}}))>0$.
There is no information of the \emph{Puk\'{a}nszky invariant} yet.
So assume that $Puk(A)=\{n_{i}: n_{i}\in \mathbb{N}_{\infty}, i\in
I\}$, where the indexing set $I$ could be finite or countable. Let
$$E_{n_{i}}=\{(s,t)\in \Delta([0,1])^{c}: m_{[0,1]}(s,t)=n_{i}\},$$
where $m_{[0,1]}$ denotes the multiplicity function of the direct
integral decomposition of $L^{2}(\mathcal{M})$ over $[0,1]\times
[0,1]$ with respect to the measure $\eta$. Then for each $i\in I$ it
is well known that $E_{n_{i}}$ are $\eta$-measurable sets. Also
\begin{align}
\nonumber &\int_{E_{n_{i}}}^{\oplus}\C_{s,t}^{n_{i}}d\eta(s,t)\cong
L^{2}(E_{n_{i}},\eta_{\mid E_{n_{i}}})\otimes \C^{n_{i}}\text{ where
}\C_{s,t}^{n_{i}}=\C^{n_{i}}, \text{ and }\\
\nonumber &\underset{{i\in I}}\oplus L^{2}(E_{n_{i}},\eta_{\mid
E_{n_{i}}})\otimes \C^{n_{i}}\cong L^{2}(\mathcal{M})\ominus
L^{2}(A).
\end{align}
In the above equation $\C^{\infty}$ stands for $l^{2}(\mathbb{N})$.
Fix orthonormal bases $\{e_{j}^{(n_{i})}\}_{1\leq j\leq n_{i}}$ of
$\C^{n_{i}}$ for all $i\in I$.\\
\noindent Then $$ \underset{i\in I}\sum \chi_{\Gamma({h_{Y}}_{\mid
F^{Y}})\cap E_{n_{i}}}\otimes e^{(n_{i})}_{1}$$ where $\chi$ denotes
the indicator function, 
can be identified with a vector $\zeta_{Y}\in
(1-e_{A})(L^{2}(\mathcal{M}))$ such that
\begin{align}\label{unitary_find}
A\zeta_{Y}A=A\zeta_{Y}=\zeta_{Y} A.
\end{align}
Eq. \eqref{unitary_find} is easy to check, in fact one only uses
that fact that $h_{Y}$ is locally one to one and onto. That
$\zeta_{Y}\neq 0$ is due to the fact $\eta(\Gamma({h_{Y}}_{\mid
F^{Y}}))>0$. Then from Theorem
\ref{identify_operators_in_q_Hilbert_normalisers}, it follows that
$\zeta_{Y}=v_{Y}\rho_{Y}$ where
$\rho_{Y}=(\zeta_{Y}^{*}\zeta_{Y})^{\frac{1}{2}}\in L^{2}(A)^{+}$
and $v_{Y}\in \mathcal{GN}(A)$. Clearly, $v_{Y}\not\in A$, as
otherwise $\overline{A\zeta_{Y}A}^{\norm{\cdot}_{2}}\subseteq
L^{2}(A)$ would become the direct integral of complex numbers over
some subset of the diagonal with respect to the measure
$\Delta_{*}\lambda$, where $\Delta:[0,1]\mapsto [0,1]\times
[0,1]$ is the map $\Delta(x)=(x,x)$.\\ 
Thus $\zeta_{Y}\in L^{2}(N(A)^{\prime\prime})$ and hence
$\overline{A\zeta_{Y}A}^{\norm{\cdot}_{2}}\subseteq
L^{2}(N(A)^{\prime\prime})$. Clearly,
\begin{align}\label{normaliser_decomposition}
\overline{A\zeta_{Y}A}^{\norm{\cdot}_{2}}\cong\int_{\Gamma({h_{Y}}_{\mid
F^{Y}})}^{\oplus}\C_{s,t}d\eta(s,t),\text{ where }\C_{s,t}=\C.
\end{align}
So $\overline{A\zeta_{Y}A}^{\norm{\cdot}_{2}}\perp L^{2}(A)$ and
$\overline{A\zeta_{Y}A}^{\norm{\cdot}_{2}}\in C_{d}(A)$. Since
$\overline{A\zeta_{Y}A}^{\norm{\cdot}_{2}}\subseteq
L^{2}(N(A)^{\prime\prime})$ so $\eta(\Gamma({h_{Y}}_{\mid
F^{Y}})\cap E_{n_{i}})=0$ if $n_{i}\geq 2$ from a result of Popa
\cite{MR815434}. Thus $1\in Puk(A)$.
%
%
%
\end{proof}


\indent Each partial isometry $0\neq v\in \mathcal{GN}(A)$
implements a measure preserving local isomorphism
$T:([0,1],\lambda)\mapsto ([0,1],\lambda)$ such that $vav^{*}=a\circ
T^{-1}$ for all $a\in A$. With abuse of notation we will write
$v=T$. Then $\Gamma(v)=\{(T(t),t): t\in Dom(T)\}$, $Dom (T)$
denoting the domain of $T$.

\begin{lem}\label{locally_decompose}
Let $\eta_{a}\neq 0$. Let $Y\subseteq (\Delta[0,1])^{c}$ be a
$\eta$-measurable set of strictly positive $\eta_{a}$-measure. Then
there is a nonzero partial isometry $v\in \mathcal{GN}(A)$ such that
$\Gamma(v)\subseteq Y$.
\end{lem}
\begin{proof}
By Lemma \ref{find_partial_isometry}, there exists a subset $F^{Y}$
of $[0,1]$ such that $\lambda(F^{Y})>0$, a invertible map $h_{Y}:
F^{Y}\mapsto h_{Y}(F^{Y})$ and a nonzero vector $\zeta_{Y}\in
L^{2}(N(A)^{\prime\prime})\ominus L^{2}(A)$ such that
$\zeta_{Y}=v_{Y}\rho_{Y}$ with $v_{Y}\in \mathcal{GN}(A)$,
$\rho_{Y}\in L^{2}(A)^{+}$ and satisfying property $(ii)$, $(iii)$,
$(iv)$ of Lemma \ref{find_partial_isometry}.\\
Let $\eta_{\zeta_{Y}}$, $\eta_{v_{Y}}$ be the measures on
$[0,1]\times [0,1]$ defined in Eq. \eqref{measure_from_kappa}. Let
$q_{Y}=v_{Y}v_{Y}^{*}\in A$. With abuse of notation we will regard
$q_{Y}$ as a measurable subset of $[0,1]$ as well. We claim that,
$\eta_{\zeta_{Y}}\ll \eta_{v_{Y}}\ll \eta_{\zeta_{Y}}$. Indeed for
$a,b\in C[0,1]$,
\begin{align}
\nonumber
\int_{[0,1]\times[0,1]}a(s)b(t)d\eta_{\zeta_{Y}}(s,t)&=\int_{\Gamma(h_{Y})}a(s)b(t)d\eta_{\zeta_{Y}}(s,t)\\
\nonumber &=\tau(\rho_{Y}^{*}v_{Y}^{*}a v_{Y}\rho_{Y} b)\\
\nonumber &=\tau(\rho_{Y}^{*}v_{Y}^{*}a v_{Y}b\rho_{Y}) \text{ (as }\rho_{Y}b=b\rho_{Y})\\
\nonumber &=\tau(v_{Y}^{*}a v_{Y}b\rho_{Y}\rho_{Y}^{*})\\
\nonumber &=\tau(v_{Y}^{*}a v_{Y}b\rho_{Y}^{*}\rho_{Y})\\
\nonumber &=\tau(v_{Y}^{*}a v_{Y}\rho_{Y}^{*}\rho_{Y}b)\\
\nonumber&=\int_{q_{Y}}a(v_{Y}(t))b(t)\abs{\rho_{Y}(t)}^{2}d\lambda(t)\\
\nonumber&=\int_{\Gamma(v_{Y})}a(s)b(t)\abs{\rho_{Y}(t)}^{2}d\eta_{v_{Y}}(s,t)\\
\nonumber&=\int_{[0,1]\times
[0,1]}a(s)b(t)\abs{\rho_{Y}(t)}^{2}d\eta_{v_{Y}}(s,t).
\end{align}
In the above string of equalities we have used the facts that $\tau$
extends to a trace like functional on $L^{1}(A)$ and the left and
right actions of $A$ on $L^{2}(A)$, $L^{1}(A)$ coincides. Using Thm.
\ref{identify_operators_in_q_Hilbert_normalisers}, by standard
arguments it follows that $\eta_{\zeta_{Y}}\ll \eta_{v_{Y}}\ll
\eta_{\zeta_{Y}}$. Thus the result follows with $v=v_{Y}$.
\end{proof}

Suppose $\{v_{j}\}_{j\in J}$ is a family of partial isometries in
$\mathcal{GN}(A)$ such that $Av_{j}\perp Av_{j^{\prime}}$ whenever
$j\neq j^{\prime}$. Denote by \cite{MR815434}
\begin{align}
\nonumber \sum_{j\in J}Av_{j}=\left\{x\in \mathcal{M}: x=\sum_{j\in
J}a_{j}v_{j}, \text{ for }a_{j}\in A\text{ with }\sum_{j\in
J}\norm{a_{j}v_{j}}_{2}^{2}<\infty\right\}.
\end{align}

\begin{thm}$($Compare Cor. $2.5$
\cite{MR815434}$)$\label{decompose_vN_normaliser} Let $\eta_{a}\neq
0$. Then $A\subsetneq N(A)^{\prime\prime}$. Moreover, there is a
sequence $\{v_{n}\}_{n=0}^{\infty}\subset \mathcal{GN}(A)$ of
nonzero partial isometries $($with possibility that the sequence
could be finite$)$ with $v_{0}=1$ such that,
\begin{align}
\nonumber &(i)\text{  }\Gamma(v_{n})\cap \Gamma(v_{m})=\emptyset \text{ for }n\neq m,\\
\nonumber &(ii)\text{ } \eta_{a}([0,1]\times [0,1])=\sum_{n=1}^{\infty}\eta_{a}(\Gamma(v_{n})),\\
\nonumber &(iii)\nonumber\oplus_{n=0}^{\infty}\overline{Av_{n}}^{\norm{\cdot}_{2}}\cong\int_{\cup_{n=0}^{\infty}\Gamma(v_{n})}^{\oplus}\C_{s,t}d(\eta_{a}+\Delta_{*}\lambda)(s,t)\cong L^{2}(N(A)^{\prime\prime}),\\
\nonumber &\indent\indent(\text{where } \C_{s,t}=\C \text{ and } \Delta:[0,1]\mapsto [0,1]\times [0,1] \text{ by }\Delta(x)=(x,x))\\
\nonumber &(iv)\text{
}N(A)^{\prime\prime}=\sum_{n=0}^{\infty}Av_{n},
\end{align}
and $\A$ restricted to
$\oplus_{n=0}^{\infty}\overline{Av_{n}}^{\norm{\cdot}_{2}}$ is
diagonalizable with respect to the decomposition in $(iii)$.
\end{thm}
\begin{proof}
First of all assuming that $(i)$ in the statement is true it follows
that $Av_{n}\perp Av_{m}$ whenever $n\neq m$. Indeed,
$\overline{Av_{n}}^{\norm{\cdot}_{2}}\subseteq
L^{2}(N(A)^{\prime\prime})$. Now $\A$ restricted to
$\overline{Av_{n}}^{\norm{\cdot}_{2}}$ is an abelian algebra with a
cyclic vector, so it is maximal abelian. The projection $e_{Av_{n}}$
onto $\overline{Av_{n}}^{\norm{\cdot}_{2}}$ is in $\A$. So
$\overline{Av_{n}}^{\norm{\cdot}_{2}}$ is the direct integral of
complex numbers over a subset $X_{n}$ of $[0,1]\times [0,1]$ with
respect to the measure $\eta$ and $\A$ restricted to
$\overline{Av_{n}}^{\norm{\cdot}_{2}}$ is diagonalizable with
respect to this decomposition. But $\eta(X_{n}\Delta
\Gamma(v_{n}))=0$. Again $\eta_{n.a}(\Gamma(v_{n}))=0$. So the
direct integral as stated above is actually with respect to the
measure $\eta_{a}+\Delta_{*}\lambda$. The graphs being disjoint for
$n\neq m$ forces the orthogonality of $Av_{n}$ and $Av_{m}$ whenever
$n\neq m$. The
sum in $(iii)$ therefore makes sense.\\
Using Lemma \ref{locally_decompose}, choose a maximal family
$\{v_{\alpha}\}_{\alpha\in \Lambda}\subset \mathcal{GN}(A)$, for
some indexing set $\Lambda$, such that $\Gamma(v_{\alpha})\subset
\Delta([0,1])^{c}$ for all $\alpha\in \Lambda$ and
$\Gamma(v_{\alpha})\cap \Gamma(v_{\beta})=\emptyset$ whenever
$\alpha\neq \beta$. Since $Av_{\alpha}\perp Av_{\beta}$ whenever
$\alpha\neq \beta$ $($by similar argument as above$)$ so the
indexing set must be countable by separability assumption of
$L^{2}(\mathcal{M})$. So we index this maximal family by
$\{v_{n}\}_{n=1}^{\infty}$. Let $v_{0}=1$. So $(i)$ follows by
construction.\\
If $\eta_{a}([0,1]\times
[0,1])>\sum_{n=1}^{\infty}\eta_{a}(\Gamma(v_{n}))$ then
$S_{a}\setminus \cup_{n=1}^{\infty}\Gamma(v_{n})$ is a set of
strictly positive $\eta_{a}$ measure. A further application of Lemma
\ref{locally_decompose} violates the maximality of
$\{v_{n}\}_{n=1}^{\infty}$. This proves $(ii)$.\\
By the argument of the first paragraph and Lemma 5.7
\cite{MR2261688},
\begin{align}\label{normaliser_maximal_decompose}
\oplus_{n=1}^{\infty}\overline{Av_{n}}^{\norm{\cdot}_{2}}\cong
\int_{\cup_{n=1}^{\infty}\Gamma(v_{n})}^{\oplus}\C_{s,t}d\eta_{a}(s,t)\subseteq
L^{2}(N(A)^{\prime\prime})\ominus L^{2}(A)
\end{align}
and $\A$ restricted to
$\oplus_{n=1}^{\infty}\overline{Av_{n}}^{\norm{\cdot}_{2}}$ is
diagonalizable with respect to the decomposition in Eq.
\eqref{normaliser_maximal_decompose}.\\
If $0\neq\zeta=\zeta^{*}\in L^{2}(N(A)^{\prime\prime})\ominus
L^{2}(A)$ is such that $\zeta\perp Av_{n}$ for all $n\geq 1$ then
$A\zeta A \perp Av_{n}$ for all $n\geq 0$. By arguments similar to
the first paragraph, $\overline{A\zeta A}^{\norm{\cdot}_{2}}$ is the
direct integral over a $\eta$-measurable set $X_{\zeta}$, of complex
numbers with respect to the measure $\eta$ and $\A$ restricted to
$\overline{A\zeta A}^{\norm{\cdot}_{2}}$ is diagonalizable
respecting this decomposition. If $\zeta$ as a $L^{2}$ function
stays nonzero on a set of positive $\Delta_{*}\lambda$-measure then
$\zeta$ cannot be perpendicular to $L^{2}(A)$. By Prop.
\ref{finite_atomic_measures}, $\overline{A\zeta
A}^{\norm{\cdot}_{2}}\in C_{d}(A)$ and hence by Theorem
\ref{fibre_invariance_theorem} and Lemma 5.7 \cite{MR2261688}, we
can assume $X_{\zeta}\subset S_{a}\setminus \Delta([0,1])$. Since
$\zeta\neq 0$ so $\eta(X_{\zeta})=\eta_{a}(X_{\zeta})>0$. Since
$\eta_{a}$ is concentrated on $\cup_{n=1}^{\infty}\Gamma (v_{n})$,
so $X_{\zeta}\cap \Gamma(v_{n})$ has strictly positive $\eta_{a}$
and hence $\eta$ measure for some $n\geq 1$. Note that
$e_{N(A)^{\prime\prime}}\in \A$ and $\A
e_{N(A)^{\prime\prime}}=\A^{\prime}e_{N(A)^{\prime\prime}}$ from
\cite{MR815434}. On the other hand, by Lemma 5.7 \cite{MR2261688},
$L^{2}(N(A)^{\prime\prime})\ominus L^{2}(A)$ will be expressed as a
direct integral over some subset of $[0,1]\times [0,1]$ with respect
to $\eta$, with multiplicity strictly bigger than $1$ on a set of
positive $\eta$-measure. This contradicts $\A
e_{N(A)^{\prime\prime}}$ is maximal abelian. Thus
\begin{align}
\nonumber\oplus_{n=0}^{\infty}\overline{Av_{n}}^{\norm{\cdot}_{2}}\cong\int_{\cup_{n=0}^{\infty}\Gamma(v_{n})}^{\oplus}\C_{s,t}d(\eta_{a}+\Delta_{*}\lambda)(s,t)\cong
L^{2}(N(A)^{\prime\prime}),
\end{align}
with associated statements about diagonalizability of $\A$. Finally
\begin{align}
\nonumber\sum_{n=0}^{\infty}Av_{n}&=\overline{\sum_{n=0}^{\infty}Av_{n}}^{\norm{\cdot}_{2}}\cap\mathcal{M}
=\left(\oplus_{n=0}^{\infty}\overline{Av_{n}}^{\norm{\cdot}_{2}}\right)\cap\mathcal{M}
=L^{2}(N(A)^{\prime\prime})\cap \mathcal{M} =N(A)^{\prime\prime}.
\end{align}
\end{proof}

\begin{rem}
Thm. \ref{decompose_vN_normaliser} generalizes Cor. 2.5 of
\cite{MR815434}. In general we cannot hope to find unitaries as was
the case in Cor. 2.5 \cite{MR815434}. The situation in Cor. $2.5$ of
\cite{MR815434} was completely different, where the assumption was
that, the masa is Cartan. Assuming the masa is Cartan, forces the
disintegration of the measure $\eta_{a}$ to have at least one atom
off the diagonal in almost every fibre. Such an assumption cannot be
made for a general masa. For example consider the following
situation. Let $C\subset \mathcal{R}$ be a Cartan masa and let
$S\subset \mathcal{R}$ be a singular masa, where $\mathcal{R}$
denotes the hyperfinite $\rm{II}_{1}$ factor. Then $C\oplus S\subset
\mathcal{R}\oplus \mathcal{R}$ is a masa, where the trace on
$\mathcal{R}\oplus \mathcal{R}$ is
$\frac{1}{2}\tau_{\mathcal{R}}\oplus \frac{1}{2}\tau_{\mathcal{R}}$,
$\tau_{\mathcal{R}}$ denoting the unique, normal, faithful tracial
state of $\mathcal{R}$. Then $C\oplus S\subset (\mathcal{R}\oplus
\mathcal{R})*\mathcal{R}\cong L(\mathbb{F}_{2})$ $($from
\cite{MR1201693}$)$ is a masa for which such an assumption will fail
from Prop. 5.10 \cite{MR2261688}.
\end{rem}

We will now present the proof of Thm \ref{classification_theorem}.
\begin{proof}[Proof of \ref{classification_theorem}.]
Case $(i)$. 
The necessary and sufficient condition for Cartan masas follows
directly from Thm. \ref{decompose_vN_normaliser}.\\
Case $(ii)$. The result for singular masas also follows from Thm.
\ref{decompose_vN_normaliser}.\\
Case $(iii)$. Let A $\varsubsetneq N(A)^{\prime\prime}\varsubsetneq
\mathcal{M}$. If $\eta_{a}=0$ then, by conclusion of $(ii)$, $A$
would become singular. Therefore $\eta_{a}\neq 0$. If $\eta_{n.a}=0$
then by conclusion of part $(i)$, $A$ would be Cartan. Therefore
$\eta_{n.a}\neq 0$ as well.\\
Conversely, if $\eta_{n.a}\neq 0$ and $\eta_{a}\neq 0$, then by
Theorem \ref{decompose_vN_normaliser}, $A\varsubsetneq
N(A)^{\prime\prime}\varsubsetneq \mathcal{M}$.\\
Case $(iv)$. First assume that $N(A)^{\prime\prime}$ is a factor.
From Thm. \ref{decompose_vN_normaliser} it follows that
$$L^{2}(N(A)^{\prime\prime})\cong\int_{[0,1]\times
[0,1]}^{\oplus}\C_{s,t}d(\eta_{a}+\tilde{\Delta}_{*}\lambda)(s,t),\text{
where }\C_{s,t}=\C$$ and $\A$ restricted to this subspace is
diagonalizable, where $\tilde{\Delta}: [0,1]\mapsto [0,1]\times
[0,1]$ is defined by $\tilde{\Delta}(x)=(x,x)$. Therefore
$[\eta_{a}+\tilde{\Delta}_{*}\lambda]$ is the
\emph{left-right-measure} for the inclusion $A\subset
N(A)^{\prime\prime}$. If the closed support of $\eta_{a}$ is
strictly contained in $[0,1]\times [0,1]$ then, there must be a open
set $U\subset [0,1]\times [0,1]$ such that $\eta_{a}(U)=0$. It
follows that the map $a\otimes b\mapsto
aJ_{N(A)^{\prime\prime}}b^{*}J_{N(A)^{\prime\prime}}$, where $a,b\in
C([0,1])$ was not an injection. On the other hand, as
$N(A)^{\prime\prime}$ is a factor the above map must be an injection
by Prop. \ref{sakai's_theorem_for_injectivity}. This
contradiction proves that the closed support of $\eta_{a}$ is $[0,1]\times [0,1]$.\\
Conversely assume $N(A)^{\prime\prime}$ has a nontrivial center. Let
$p\in \textbf{Z}(N(A)^{\prime\prime})$ be a projection which is
different from $0$ and $1$. Then
$$N(A)^{\prime\prime}=N(A)^{\prime\prime}p\oplus N(A)^{\prime\prime}(1-p).$$
So $p\in A^{\prime}\cap N(A)^{\prime\prime}$ and hence $p\in A$. So
\begin{align}\label{direct_sum_decompose}
A=Ap\oplus A(1-p).
\end{align}
It follows that are exists $\lambda$-measurable sets
$F_{1},F_{2}\subset [0,1]$ such that, $F_{1}\cup F_{2}=[0,1]$,
$F_{1}\cap F_{2}=\emptyset$, $C(F_{1})$, $C(F_{2})$ are \emph{w.o.t}
dense unital subalgebras of $Ap$ and $A(1-p)$ respectively and
$C(F_{1})\oplus C(F_{2})$ is \emph{w.o.t} dense in $A$. With respect
to the Eq. \eqref{direct_sum_decompose} let $a=a_{1}\oplus a_{2}$
and $b=b_{1}\oplus b_{2}$ be the decompositions of $a,b\in
C([0,1])$,
where $a_{i},b_{i}\in C(F_{i})$ for $i=1,2$. 
For $\zeta\in L^{2}(N(A)^{\prime\prime})$ one has an analogous
decomposition $\zeta=\zeta_{1}\oplus\zeta_{2}$ with
$\zeta_{1}=p\zeta$ and $\zeta_{2}=(1-p)\zeta$. The equation
\begin{align}
\nonumber \langle a\zeta b,\zeta\rangle &= \langle (a_{1}\oplus
a_{2})(\zeta_{1}\oplus\zeta_{2})(b_{1}\oplus
b_{2}),(\zeta_{1}\oplus\zeta_{2})\rangle =\langle a_{1}\zeta_{1}
b_{1},\zeta_{1}\rangle+\langle a_{2}\zeta_{2} b_{2},\zeta_{2}\rangle
\end{align}
shows that the \emph{left-right-measure} for the inclusion $A\subset
N(A)^{\prime\prime}$ will be concentrated on $F_{1}\times F_{1} \cup
F_{2}\times F_{2}$. It follows that closed support of $\eta_{a}$ is
strictly contained in $[0,1]\times [0,1]$. This completes the proof.
\end{proof}

\begin{rem}\label{not_injection}
The proof of Case $(iv)$ actually shows that if $A\subset
\mathcal{N}$ is a masa where $\mathcal{N}$ is a finite type
$\rm{II}$ algebra with a nontrivial center then for any choice of
compact Hausdorff space $X$ such that $C(X)$ is \emph{w.o.t} dense,
unital and norm separable subalgebra of $A$, the map
$$\overset{n}{\underset{i=1}\sum}a_{i}\otimes b_{i}\mapsto
\overset{n}{\underset{i=1}\sum}a_{i}Jb_{i}^{*}J$$ from
$C(X)\otimes_{alg} C(X)\mapsto \textbf{B}(L^{2}(\mathcal{N}))$ is
not an injection for any choice of trace.
\end{rem}

The following results about masas that were proved by experts in
different ways are just easy consequences of the \emph{measurable
selection principle} as we have described in this section.

\begin{cor}\label{Cartan_in_intermidiate}
If $A\subset \mathcal{M}$ is a Cartan masa then $A\subset B$ is a
Cartan masa for all von Neumann subalgebra $A\varsubsetneq
B\varsubsetneq \mathcal{M}$.
\end{cor}
\begin{proof}
By Lemma 5.7 of \cite{MR2261688} the \emph{left-right-measure} of
the inclusion $A\subset \mathcal{M}$ is
$[\eta_{B}+\eta_{B^{\perp}}]$ where $[\eta_{B}]$ is the
\emph{left-right-measure} of the inclusion $A\subset B$ and
$\eta_{B}\perp \eta_{B^{\perp}}$. It follows that $\eta_{B}$ has
atomic disintegration along both axes. The result is then immediate
from Thm. \ref{classification_theorem} and Thm
\ref{decompose_vN_normaliser}.
\end{proof}

\begin{cor}
Let $A\subset \mathcal{M}$ be a masa and let $Q$ be a finite von
Neumann algebra such that $dim(Q)\geq 2$. Then
$N_{\mathcal{M}*Q}(A)=N_{\mathcal{M}}(A)$.
\end{cor}
\begin{proof}
In this proof we consider \emph{left-right-measures} restricted to
the off diagonal. First of all it well known that $A\subset
\mathcal{M}*Q$ is a masa. Let $[\eta_{\mathcal{M}}]$ denote the
\emph{left-right-measure} of the inclusion $A\subset \mathcal{M}$.
Write $\eta_{\mathcal{M}}=\eta_{1}+\eta_{2}$ where $\eta_{1}\ll
\lambda\otimes \lambda$ and $\eta_{2}\perp \lambda\otimes \lambda$.
Using Prop. $5.10$ and Lemma $5.7$ \cite{MR2261688} it follows that
the \emph{left-right-measure} $[\eta_{\mathcal{M}*Q}]$ of the
inclusion $A\subset \mathcal{M}*Q$ is given by
\begin{equation}
\nonumber \eta_{\mathcal{M}*Q}= \begin{cases}
\eta_{\mathcal{M}}+\lambda\otimes\lambda &\text{ if } \eta_{1}= 0, \\
\eta_{2}+ \lambda\otimes\lambda &\text{ if } \eta_{1}\neq 0.
\end{cases}
\end{equation}
The rest is obvious from Thm. \ref{classification_theorem} and Thm.
\ref{decompose_vN_normaliser}.
\end{proof}

\begin{cor}
Let $A\subset \mathcal{M}$ be a Cartan masa and let $A\subset
B\subsetneq \mathcal{M}$ be an intermediate subalgebra. Then there
is a $v\in \mathcal{GN}(A)$ such that $v\perp B$.
\end{cor}
\begin{proof}
By Lemma 5.7 \cite{MR2261688}, the \emph{left-right-measure} of the
inclusion $A\subset \mathcal{M}$ is $[\eta_{B}+\eta_{B^{\perp}}]$
where $\eta_{B}\perp\eta_{B^{\perp}}$ and $[\eta_{B}]$ is the
\emph{left-right-measure} of the inclusion $A\subset B$. Note
$\eta_{B^{\perp}}\neq 0$. Apply Lemma \ref{locally_decompose}.
\end{proof}


%

We prove the next theorem in the context of $\rm{II}_{1}$ factors.
But it can be easily generalized to finite von Neumann algebras. Let
$\mathcal{M}_{i}$, $i=1,2$ be separably acting $\rm{II}_{1}$ factors
with normal, faithful tracial states $\tau_{i}$ respectively. Let
$\mathcal{M}_{i}$ act on $L^{2}(\mathcal{M}_{i},\tau_{i})$ by left
multiplication. Let $A\subset \mathcal{M}_{1}$ and $B\subset
\mathcal{M}_{2}$ be masas. Fix compact Polish spaces $X,Y$ such that
$C(X)\subset A$ and $C(Y)\subset B$ are unital, norm separable and
\emph{w.o.t} dense. Let $\nu_{X}$ and $\nu_{Y}$ denote the tracial
measures for $A,B$ respectively, which will be assumed to be
complete. Let the \emph{left-right-measure} of $A$ on $X\times X$ be
$[\sigma_{1}]$ and that of $B$ on $Y\times Y$ be $[\sigma_{2}]$.
Here we are allowing the diagonals, i.e. we are assuming
${\sigma_{1}}_{\mid\Delta(X)}=(\tilde{\Delta}_{X})_{*}\nu_{X}$ and
${\sigma_{2}}_{\mid\Delta(Y)}=(\tilde{\Delta}_{Y})_{*}\nu_{Y}$ where
$\tilde{\Delta}_{X}:X\mapsto X\times X$ by
$\tilde{\Delta}_{X}(x)=(x,x)$ and $\tilde{\Delta}_{Y}:Y\mapsto
Y\times Y$ by $\tilde{\Delta}_{Y}(y)=(y,y)$.\\
\indent By Tomita's theorem on commutants $A\overline{\otimes} B$ is
a masa in $\mathcal{M}_{1}\overline{\otimes}\mathcal{M}_{2}$.
$X\times Y$ is compact and Polish, and $C(X\times Y)$ is unital,
norm separable and \emph{w.o.t} dense in $A\overline{\otimes} B$.
The standard Hilbert space and the Tomita's involution operator for
$\mathcal{M}_{1}\overline{\otimes}\mathcal{M}_{2}$ is
$L^{2}(\mathcal{M}_{1},\tau_{1})\otimes
L^{2}(\mathcal{M}_{2},\tau_{2})$ and $J_{\mathcal{M}_{1}}\otimes
J_{\mathcal{M}_{2}}$ respectively. The tracial measure for
$A\overline{\otimes} B$ on $X\times Y$ is clearly $\nu_{X}\otimes
\nu_{Y}$. With this as set up we formulate the next theorem which
appeared in \cite{MR999997}. The same proof actually generalizes to
infinite tensor products.

\begin{thm}$($Chifan's Normaliser Formula$)$\label{Chifan_theorem}
Let $A\subset \mathcal{M}_{1}$ and $B\subset \mathcal{M}_{2}$ be
masas in separably acting $\rm{II}_{1}$ factors $\mathcal{M}_{1}$
and $\mathcal{M}_{2}$. Then
$$N(A\overline{\otimes} B)^{\prime\prime}=N(A)^{\prime\prime}\overline{\otimes} N(B)^{\prime\prime}.$$
\end{thm}
\begin{proof}
Fix $\sigma_{1}$ and $\sigma_{2}$ from the aforesaid class of
\emph{left-right-measures}. The \emph{left-right-measure} of
$A\overline{\otimes} B$ on $(X\times Y)\times (X\times Y)$ which is
denoted by $[\beta]$ is given by
$$d\beta(s_{X},s_{Y},t_{X},t_{Y})=d\sigma_{1}(s_{X},t_{X})d\sigma_{2}(s_{Y},t_{Y})
$$
from Prop. $5.2$ \cite{MR2261688}. Here $s$ is the variable running
along the first coordinate $($horizontal direction$)$ and $t$ along
the second coordinate $($vertical direction$)$. Then from Lemma
\ref{fibre_tensor_formulae} it follows that the disintegration of
$\beta$ along the $t$ variable $($vertical direction$)$ is given by
$$\beta_{t_{X\times Y}}={\sigma_{1}}_{t_{X}}\otimes {\sigma_{2}}_{t_{Y}}, \text{ } (t_{X},t_{Y}) \text{ - a.e }\nu_{X}\otimes \nu_{Y}, \text{ where }t_{X\times Y}=(t_{X},t_{Y}).$$
For fixed $t_{X\times Y}=(t_{X},t_{Y})\in X\times Y$ the measure
$\beta_{t_{X\times Y}}$ has an atom at the point
$(s_{X},s_{Y},t_{X},t_{Y})$ if and only if ${\sigma_{1}}_{t_{X}}$
has an atom at $(s_{X},t_{X})$ and ${\sigma_{2}}_{t_{Y}}$ has an
atom at $(s_{Y},t_{Y})$. Therefore $$((X\times Y)\times (X\times
Y))_{a}=S_{2,3}((X\times X)_{a}\times (Y\times Y)_{a}),$$ where
$S_{2,3}$ denotes the permutation $(2,3)$ on four symbols $($see
Prop. \ref{atomic_part_measurable_set}$)$. Therefore,
$$\beta_{\mid ((X\times Y)\times (X\times
Y))_{a}}={\sigma_{1}}_{\mid(X\times X)_{a}}\otimes
{\sigma_{2}}_{\mid(Y\times Y)_{a}}.$$ Hence denoting
$\C_{s_{X},s_{Y},t_{X},t_{Y}}=\C$,
$\C_{s_{X},t_{X}}=\C=\C_{s_{Y},t_{Y}}$ we have
\begin{align}
\nonumber L^{2}(N(A\overline{\otimes}B)^{\prime\prime})&\cong\int_{((X\times Y)\times (X\times Y))_{a}}^{\oplus}\C_{s_{X},s_{Y},t_{X},t_{Y}}d\beta(s_{X},s_{Y},t_{X},t_{Y})\\
\nonumber &\cong\int_{(X\times X)_{a}}^{\oplus}\C_{s_{X},t_{X}}d\sigma_{1}(s_{X},t_{X})\otimes\int_{(Y\times Y)_{a}}^{\oplus}\C_{s_{Y},t_{Y}}d\sigma_{2}(s_{Y},t_{Y})\\
\nonumber &\cong L^{2}(N(A)^{\prime\prime})\otimes
L^{2}(N(B)^{\prime\prime}) \text{ from Thm.
}\ref{decompose_vN_normaliser}.
\end{align}
Since the containment $N(A)^{\prime\prime}\overline{\otimes}
N(B)^{\prime\prime}\subseteq N(A\overline{\otimes}
B)^{\prime\prime}$ is obvious we are done.
\end{proof}

\begin{cor}\cite{MR999995}
Let $A\subset \mathcal{M}_{1}$ and $B\subset \mathcal{M}_{2}$ be
singular masas in separably acting $\rm{II}_{1}$ factors
$\mathcal{M}_{1}$ and $\mathcal{M}_{2}$. Then $A\overline{\otimes}
B$ is singular in $\mathcal{M}_{1}\overline{\otimes}
\mathcal{M}_{2}$.
\end{cor}

\section{Asymptotic Homomorphism and Measure Theory }

The equivalence of WAHP and singularity is a nontrivial theorem
\cite{MR999995}. In this section we will give a direct proof of the
equivalence of WAHP and singularity by using measure theoretic
tools. We will also present partial results about AHP. In order to
do so we will first have to relate certain norms to the
\emph{left-right-measure}. The measure theoretic tools described in
this section will be used in a future paper
for explicit calculation of \emph{left-right-measures}.\\
\indent Let $A\subset \mathcal{M}$ be a masa. Let $\lambda$ denote
the Lebesgue measure on $[0,1]$ so that $A\cong
L^{\infty}([0,1],\lambda)$. Then $\lambda$ is the tracial measure.
Let $[\eta]$ denote the \emph{left-right-measure} for $A$. We assume
that $\eta$ is a probability measure on $[0,1]\times [0,1]$ and
$\eta(\Delta([0,1]))=0$. Let $B[0,1]$ denote the collection of
all bounded measurable functions on $[0,1]$.\\
\noindent \textbf{Notation:} The disintegrated measures are usually
written with a subscript $t\mapsto \eta_{t}$ in the literature. But
in this section we will use the superscript notation $t\mapsto
\eta^{t}$ to denote them. The $(\pi_{1},\lambda)$ disintegration of
measures will be indexed by the variable $t$ and the $(\pi_{2},\lambda)$
disintegration will be indexed by the variable $s$.\\
\indent In all the following results that uses disintegration of
measures we will only state or prove the result with respect to the
$(\pi_{1},\lambda)$ disintegration. Statements about the
$(\pi_{2},\lambda)$ disintegration are analogous.

\begin{lem}\label{identify_disintegrated_measure}
Let $x\in \mathcal{M}$ be such that $\mathbb{E}_{A}(x)=0$. Let
$\eta_{x}$ denote the measure on $[0,1]\times [0,1]$ defined in Eq.
\eqref{measure_from_kappa}. Then $\eta_{x}$ admits
$(\pi_{i},\lambda)$ disintegrations $[0,1]\ni t\mapsto\eta_{x}^{t}$
and $[0,1]\ni s\mapsto\eta_{x}^{s}$, where $\pi_{i}$, $i=1,2$
denotes the coordinate projections. Moreover,
\begin{align}
\nonumber \eta_{x}^{t}([0,1]\times [0,1])=\mathbb{E}_{A}(xx^{*})(t),
\text{ }\lambda \text{ a.e.}
\end{align}
\end{lem}
\begin{proof}
From Lemma 5.7 of \cite{MR2261688} it follows that there is a
measure $\eta_{0}$ such that $(i)$ $\eta_{0}\perp \eta_{x}$, $(ii)$
$[\eta_{x}+\eta_{0}]$ is the \emph{left-right-measure} for $A$.
Therefore $[(\pi_{i})_{*}(\eta_{0}+\eta_{x})]=[\lambda]$ by Lemma
\ref{properties_of_eta} and hence $(\pi_{i})_{*}(\eta_{x})\ll
\lambda$ for $i=1,2$. Consequently from Thm.
\ref{Existence_of_disintegration},
$\eta_{x}$ admits $(\pi_{i},\lambda)$ disintegrations for $i=1,2$.\\
Note that
$\eta_{x}([0,1]\times[0,1])=\tau(xx^{*})=\tau(\mathbb{E}_{A}(xx^{*}))$.
From $(ii)$ of Defn. \ref{definition_of_disintegration} it follows
that $[0,1]\ni t\mapsto \eta_{x}^{t}([0,1]\times [0,1])$ is
measurable. Let $E\subseteq [0,1]$ be any Borel set. Then there
exists a sequence of functions $f_{n}\in C[0,1]$ such that $0\leq
f_{n}\leq 1$ and $f_{n}\rightarrow \chi_{E}$ pointwise. By dominated
convergence theorem we have $\eta_{x}(f_{n}\otimes 1)\rightarrow
\eta_{x}(\chi_{E}\otimes 1)$. On the other hand,
\begin{align}
\nonumber\eta_{x}(f_{n}\otimes 1)&=\langle f_{n}x,x\rangle
=\tau(f_{n}xx^{*}) =\tau(f_{n}\mathbb{E}_{A}(xx^{*}))
=\int_{0}^{1}f_{n}(t)\mathbb{E}_{A}(xx^{*})(t)d\lambda(t)\\
\nonumber &\rightarrow\int_{0}^{1}\chi_{E}(t)\mathbb{E}_{A}(xx^{*})(t)d\lambda(t)\text{, }\text{as }n\rightarrow \infty,\\
\nonumber &=\int_{E}\mathbb{E}_{A}(xx^{*})(t)d\lambda(t).
\end{align}
From Defn. \ref{definition_of_disintegration} again we have
\begin{align}
\nonumber \eta_{x}(\chi_{E}\otimes
1)&=\int_{0}^{1}\eta_{x}^{t}(\chi_{E}\otimes 1)d\lambda(t)
=\int_{E}\eta_{x}^{t}([0,1]\times [0,1])d\lambda(t).
\end{align}
Therefore for all Borel sets $E\subseteq [0,1]$ we have
\begin{align}
\nonumber\int_{E}\eta_{x}^{t}([0,1]\times
[0,1])d\lambda(t)=\int_{E}\mathbb{E}_{A}(xx^{*})(t)d\lambda(t).
\end{align}
Thus, $\eta_{x}^{t}([0,1]\times [0,1])=\mathbb{E}_{A}(xx^{*})(t)$
for $\lambda$ almost all $t$.
\end{proof}

\begin{lem}\label{integrability}
Let $x\in \mathcal{M}$ be such that $\mathbb{E}_{A}(x)=0$. Let $f\in
B[0,1]$. Then the functions $[0,1]\ni t\mapsto \eta_{x}^{t}(1\otimes
f), [0,1]\ni s\mapsto \eta_{x}^{s}(f\otimes 1)$ are in
$L^{\infty}([0,1],\lambda)$.
\end{lem}
\begin{proof}
We will only prove for the $(\pi_{1},\lambda)$ disintegration. From
Lemma \ref{identify_disintegrated_measure} we know that $\eta_{x}$
admits a $(\pi_{1},\lambda)$ disintegration. From Defn.
\ref{definition_of_disintegration} we also know that $[0,1]\ni
t\mapsto \eta_{x}^{t}(1\otimes f)$ is measurable. Now of $0\leq
t\leq 1$
\begin{align}
\nonumber \abs{\eta_{x}^{t}(1\otimes f)}\leq
\norm{f}\eta_{x}^{t}([0,1]\times [0,1]).
\end{align}
Now use Lemma \ref{identify_disintegrated_measure}.
\end{proof}

\begin{lem}\label{relate_norm_to_fourier_coeff}
Let $x\in \mathcal{M}$ be such that $\mathbb{E}_{A}(x)=0$. Let
$b,w\in B[0,1]$. Then
\begin{align}
\nonumber
\norm{\mathbb{E}_{A}(bxwx^{*})}_{2}^{2}=\int_{0}^{1}\abs{b(t)}^{2}\abs{{\eta_{x}^{t}(1\otimes
w)}}^{2}d\lambda(t).
\end{align}
\end{lem}
\begin{proof}
We have noted before that $\eta_{x}$ admits $(\pi_{i},\lambda)$
disintegrations for $i=1,2$. Secondly, as $b,w\in B[0,1]$, so
$[0,1]\ni t\mapsto b(t){\eta_{x}^{t}(1\otimes w)}$ is in
$L^{\infty}([0,1],\lambda)$ from Lemma \ref{integrability}. Now
\begin{align}
\nonumber\norm{\mathbb{E}_{A}(bxwx^{*})}_{2}^{2}&=\underset{\norm{a}_{2}\leq
1}{\underset{a\in C[0,1]}\sup}\abs{\langle
a,\mathbb{E}_{A}(bxwx^{*})\rangle }^{2}\\
\nonumber & =\underset{\norm{a}_{2}\leq 1}{\underset{a\in
C[0,1]}\sup}\abs{\tau(a\mathbb{E}_{A}(bxwx^{*}))}^{2}\\
\nonumber &=\underset{\norm{a}_{2}\leq 1}{\underset{a\in
C[0,1]}\sup}\abs{\tau(\mathbb{E}_{A}(abxwx^{*}))}^{2}\\
\nonumber &=\underset{\norm{a}_{2}\leq 1}{\underset{a\in C[0,1]}\sup}\abs{\tau(abxwx^{*})}^{2} \\
\nonumber &=\underset{\norm{a}_{2}\leq 1}{\underset{a\in C[0,1]}\sup}\abs{\int_{[0,1]\times[0,1]}a(t)b(t)w(s)d\eta_{x}(t,s)}^{2}\text{ (from Eq. }\eqref{measure_from_kappa})\\
\nonumber &=\underset{\norm{a}_{2}\leq 1}{\underset{a\in C[0,1]}\sup}\abs{\int_{0}^{1}a(t)b(t){\eta_{x}^{t}}(1\otimes w)d\lambda(t)}^{2}\text{ }(\text{from Defn. } \ref{definition_of_disintegration})\\
\nonumber &=\int_{0}^{1}\abs{b(t)}^{2}\abs{{\eta_{x}^{t}}(1\otimes
w)}^{2}d\lambda(t) \text{ }(\text{from Lemma }\ref{integrability}).
\end{align}
\end{proof}
The following facts are well known, we just record them for
completeness. For details we refer the reader to \cite{MR2186251}.
Recall that a subset $S\subseteq \mathbb{Z}$ is said to be of
\emph{full density} if
\begin{align}
\nonumber \underset{n}\lim \frac{\# (S\cap [-n,n])}{2n+1}=1.
\end{align}

\begin{defn}\label{define_mixing_measure}
A measure $\mu$ on $[0,1]$ is called \emph{mixing} (or sometimes
\emph{Rajchman}) if its Fourier coefficients
$\hat{\mu}_{n}=\int_{0}^{1}e^{2\pi int}d\mu(t)$ converge to $0$ as
$\abs{n}\rightarrow \infty$.
\end{defn}
By the Riemann-Lebesgue lemma any absolutely continuous measure is
mixing. However there are many mixing singular measures as well.
Atomic measures can never be mixing.
The next proposition justifies why non-atomic measures are called
\emph{weak} (or \emph{weakly}) \emph{mixing measures}.
\begin{prop}$($Wiener$)$\label{weakmixingmeasures}
A measure $\mu$ on $[0,1]$ is non-atomic $($diffuse$)$ if and only
if for a set $S\subseteq \mathbb{Z}$ of full density
\begin{align}
\nonumber \underset{n\in S, \abs{n}\rightarrow  \infty}\lim
\hat{\mu}_{n}=0.
\end{align}
\end{prop}


From Prop. $2.5$ and Prop. $2.19$ of \cite{MR2186251}, \emph{mixing}
and \emph{weakly mixing} are just not properties of measures, they
are in fact properties of equivalence
class of measures.\\
\indent We need the following fact from the calculus course. A
bounded sequence of complex numbers $\{a_{n}\}_{n\in \mathbb{Z}}$
converges to $0$ strongly in the sense of Ces\`{a}ro i.e.
\begin{align}\label{cesaro}
\underset{N\rightarrow \infty}\lim\text{
}\frac{1}{2N+1}\sum_{n=-N}^{N}\abs{a_{n}}=0
\end{align}
if and only if there is a set $S\subseteq \mathbb{Z}$ of full
density such that
\begin{align}\label{converge_in_density_one}
\underset{n\in S, \abs{n} \rightarrow \infty}\lim \abs{a_{n}}=0.
\end{align}

Let $x,y\in \mathcal{M}$ be such that
$\mathbb{E}_{A}(x)=\mathbb{E}_{A}(y)=0$. Let $a\in A$. Then the
following polarization identity holds:
\begin{align}\label{polarize_cond_exp}
4\text{ }\mathbb{E}_{A}(xay^{*})&=\mathbb{E}_{A}((x+y)a(x+y)^{*})-\mathbb{E}_{A}((x-y)a(x-y)^{*})\\
\nonumber &+ i\text{ }\mathbb{E}_{A}((x+iy)a(x+iy)^{*})-i\text{
}\mathbb{E}_{A}((x-iy)a(x-iy)^{*}).
\end{align}
Thus WAHP for a masa is equivalent to the following. For each finite
set $\{x_{i}\}_{i=1}^{n}\subset \mathcal{M}$ with
$\mathbb{E}_{A}(x_{i})=0$ for all $1\leq i\leq n$ and $\epsilon>0$,
there exists an unitary $u\in A$ such that
\begin{align}
\nonumber \norm{\mathbb{E}_{A}(x_{i}ux_{i}^{*})}_{2}\leq \epsilon
\text{ for all }1\leq i\leq n.
\end{align}

We will only prove the harder part of the equivalence of singularity
and WAHP.

\begin{thm}\label{measuretoWAHP}
Let $A\subset \mathcal{M}$ be a masa such that
${L^{2}(A)}^{\perp}\in C_{n.a}(A)$. Then $A$ has WAHP.
\end{thm}
\begin{proof}
Suppose to the contrary $A$ does not have WAHP. Then there is a
$\epsilon >0$ and operators $0\neq x_{i}\in \mathcal{M}$, $1\leq
i\leq n$ with $\mathbb{E}_{A}(x_{i})=0$ for all $i$, such that
\begin{align}
\nonumber \underset{u\in \mathcal{U}(A)}\inf\sum_{i=1}^{n}
\norm{\mathbb{E}_{A}(x_{i}ux_{i}^{*})}_{2}^{2}\geq \epsilon,
\end{align}
where $\mathcal{U}(A)$ denotes the unitary group of $A$. Note that
for all $1\leq i\leq n$, $\overline{Ax_{i}A}^{\norm{\cdot}_{2}}\in
C_{n.a}(A)$ by Lemma 5.7 of \cite{MR2261688} and Thm.
\ref{classification_theorem}. Equivalently, if $t\mapsto
\eta_{x_{i}}^{t}$ and $s\mapsto \eta_{x_{i}}^{s}$ denote the
$(\pi_{1},\lambda)$ and $(\pi_{2},\lambda)$ disintegrations
respectively of $\eta_{x_{i}}$, then for $\lambda$ almost all $t$,
the measure $\eta_{x_{i}}^{t}$ is completely non-atomic and similar
statements hold for $\eta_{x_{i}}^{s}$.\\
Let $v\in A$ be the Haar unitary corresponding to the function
$t\mapsto e^{2\pi it}$. Then $v$ generates $A$. Now from Lemma
\ref{relate_norm_to_fourier_coeff} we have
\begin{align}\label{bigger_than_epsilon}
\sum_{i=1}^{n}\norm{\mathbb{E}_{A}(x_{i}v^{k}x_{i}^{*})}_{2}^{2}=\int_{0}^{1}\sum_{i=1}^{n}\abs{\eta_{x_{i}}^{t}(1\otimes
v^{k})}^{2}d\lambda(t)\geq \epsilon \text{ for all }k\in \mathbb{Z}.
\end{align}
Throwing off a $\lambda$-null set $F$ we assume that for $t\in
F^{c}$ the measures $\eta_{x_{i}}^{t}$ are completely non-atomic,
finite, concentrated on $\{t\}\times [0,1]$ and
$\eta_{x_{i}}^{t}([0,1]\times
[0,1])=\mathbb{E}_{A}(x_{i}x_{i}^{*})(t)$ for all $1\leq i\leq n$
$($see Lemma \ref{identify_disintegrated_measure}$)$. Let
\begin{align}
\nonumber a_{k}(t)=\sum_{i=1}^{n}\abs{\eta_{x_{i}}^{t}(1\otimes
v^{k})}^{2}, \text{ }k\in \mathbb{Z}, t\in [0,1].
\end{align}
Then $a_{k}$ is measurable for all $k\in \mathbb{Z}$. For $k\in
\mathbb{Z}$ and $t\in F^{c}$ we have
\begin{align}
\nonumber a_{k}(t)&=\sum_{i=1}^{n}\abs{\int_{[0,1]\times
[0,1]}e^{2\pi iks}d\eta_{x_{i}}^{t}(t^{\prime},s)}^{2} \leq
\sum_{i=1}^{n}\left(\eta_{x_{i}}^{t}([0,1]\times [0,1])\right)^{2}.
\end{align}
Then by Lemma \ref{identify_disintegrated_measure}, $a_{k}(t)\leq
\sum_{i=1}^{n}\abs{\mathbb{E}_{A}(x_{i}x_{i}^{*})(t)}^{2}<\infty$,
for all $t\in F^{c}$ and for all $k\in \mathbb{Z}$. Define
\begin{align}
\nonumber s_{N}(t)&=\frac{1}{2N+1}\sum_{k=-N}^{N}a_{k}(t), N\in
\mathbb{N}.
\end{align}
Then $s_{N}$ is measurable for all $N\in \mathbb{N}$. Since
$\eta_{x_{i}}^{t}$ is completely non-atomic for all $1\leq i\leq n$
and $t\in F^{c}$ so
\begin{align}
\nonumber s_{N}(t)\rightarrow 0\text{ as }N\rightarrow \infty\text {
for all }t\in F^{c} \text{ from Eq. } \eqref{cesaro},
\eqref{converge_in_density_one} \text{ and Prop }
\ref{weakmixingmeasures}.
\end{align}
Again since $s_{N}(t)\leq\sum_{i=1}^{n}
\abs{\mathbb{E}_{A}(x_{i}x_{i}^{*})(t)}^{2}$ for $t\in F^{c}$
$($from Lemma \ref{identify_disintegrated_measure}$)$, so by
dominated convergence theorem
\begin{align}
\nonumber \int_{0}^{1}s_{N}(t)d\lambda(t)\rightarrow 0 \text{ as
}N\rightarrow \infty.
\end{align}
Therefore,
\begin{align}
\nonumber
\int_{0}^{1}s_{N}(t)d\lambda(t)&=\frac{1}{2N+1}\sum_{k=-N}^{N}\int_{0}^{1}\sum_{i=1}^{n}\abs{\eta_{x_{i}}^{t}(1\otimes
v^{k})}^{2}d\lambda(t)\\
\nonumber&=\frac{1}{2N+1}\sum_{k=-N}^{N}\left(\sum_{i=1}^{n}
\norm{\mathbb{E}_{A}(x_{i}v^{k}x_{i}^{*})}_{2}^{2}\right)\rightarrow
0\text{ as }N\rightarrow \infty.
\end{align}
Consequently from Eq. \eqref{converge_in_density_one} there is a set
$S\subseteq \mathbb{Z}$ of full density such that
\begin{align}
\nonumber \underset{k\in S,\abs{k}\rightarrow \infty}\lim
\sum_{i=1}^{n} \norm{\mathbb{E}_{A}(x_{i}v^{k}x_{i}^{*})}_{2}^{2}=0.
\end{align}
This is a contradiction to Eq. \eqref{bigger_than_epsilon}. So $A$
must have WAHP.
\end{proof}

The proof of Thm. \ref{measuretoWAHP} yields the following result.

\begin{thm}\label{avg_go_to_zero}
Let $A\subset \mathcal{M}$ be a singular masa. Then given any finite
set $\{x_{i}\}_{i=1}^{n}\subset \mathcal{M}$ with
$\mathbb{E}_{A}(x_{i})=0$ for all $i$,
\begin{align}\label{like_weak_mixing}
\frac{1}{2N+1}\sum_{k=-N}^{N}\left(\sum_{i=1}^{n}
\norm{\mathbb{E}_{A}(x_{i}v^{k}x_{i}^{*})}_{2}^{2}\right)\rightarrow
0\text{ as }N\rightarrow \infty.
\end{align}
where $v$ is a Haar unitary generator of $A$.
\end{thm}

\begin{rem}
Thus the unitary in the definition of WAHP $($Defn.
\ref{strongly_singular_and_WAHP}$)$ can always be chosen to be
$v^{k}$ where $k$ is a large integer and $v$ is a Haar unitary
generator of the masa. This strengthens the definition of WAHP. Note
that Eq. \eqref{like_weak_mixing} is very closely related to
definition of weakly mixing actions of abelian groups on finite von
Neumann algebras.
\end{rem}

The measures $\eta_{x}^{t},\eta^{t}$ are concentrated on
$\{t\}\times [0,1]$ for $\lambda$ almost all $t$. We will denote by
$\tilde{\eta}_{x}^{t},\tilde{\eta}^{t}$ the restriction of the
measures $\eta_{x}^{t}$ and $\eta^{t}$ respectively on $\{t\}\times
[0,1]$. Thus $\tilde{\eta}_{x}^{t},\tilde{\eta}^{t}$ can be regarded
as measures on $[0,1]$.

\begin{thm}\label{mixing_to_AHP}
Let $A\subset \mathcal{M}$ be a masa. Let $[\eta]$ denote the
left-right-measure for $A$. If for $\lambda$ almost all $t$ the
measures $\tilde{\eta}^{t}$ are mixing, then $A$ has $AHP$ with
respect to a Haar unitary generator of $A$.
\end{thm}
\begin{proof}
From Prop. $2.5$ of \cite{MR2186251} it follows that for $\lambda$
almost all $t$, any measure in the equivalence class
$[\tilde{\eta}^{t}]$ is mixing. In view of Eq.
\eqref{polarize_cond_exp}, it is enough to show that for all $x\in
\mathcal{M}$ with $\mathbb{E}_{A}(x)=0$,
\begin{align}
\nonumber \norm{\mathbb{E}_{A}(xv^{n}x^{*})}_{2}\rightarrow 0 \text{
as } \abs{n}\rightarrow \infty,
\end{align}
where $v\in A$ is a Haar unitary generator of $A$. Let $v\in A$
correspond to the function $s\mapsto e^{2\pi is}$. By Lemma
\ref{relate_norm_to_fourier_coeff}
\begin{align}
\nonumber
\norm{\mathbb{E}_{A}(xv^{n}x^{*})}_{2}^{2}=\int_{0}^{1}\abs{\eta_{x}^{t}(1\otimes
v^{n})}^{2}d\lambda(t).
\end{align}
From Lemma 5.7 \cite{MR2261688} we know that $\eta_{x}\ll \eta$ and
hence for $\lambda$ almost all $t$, $\eta_{x}^{t}\ll \eta^{t}$  from
Lemma
\ref{equivalence_of_measure_imply_equivalence_of_fibre_almost_everywhere}.
So $\tilde{\eta}_{x}^{t}\ll \tilde{\eta}^{t}$ for $\lambda$ almost
all $t$. Thus $\tilde{\eta}_{x}^{t}$ is mixing measure from Prop.
$2.5$ of \cite{MR2186251} for $\lambda$ almost all $t$. Also from
Lemma \ref{identify_disintegrated_measure}, the measures
$\eta_{x}^{t}$ are finite for $\lambda$ almost all $t$. Use Lemma
\ref{identify_disintegrated_measure} and apply dominated convergence
theorem to finish the proof.
\end{proof}

\appendix
\numberwithin{equation}{section}

\section{Structure of measurable functions}

Making a measurable selection as we attempted in Lemma
\ref{construction_of_Baire_functions} is not enough. One likes to
make a measurable selection so that the graph of the selection is an
automorphism graph of the masa, the automorphism being implemented
by an unitary in the factor. But this is a very delicate issue. We
are not aware of such selection theorems. We can overcome this
obstacle though. Structure theorems of continuous and measurable
functions are what comes into play.
\begin{defn}\label{condition_(N)}
Let $f:[0,1]\mapsto \mathbb{R}$ be a function and $E$ be a subset of
[0,1]. Then $f$ is said to satisfy condition $(N)$ or null condition
of Lusin relative to $E$ if $f(A)$ is a set of measure $0$ whenever
$A\subset E$ is a set of measure $0$.
\end{defn}
\indent The definition implicitly assumes that there are two
measures on $[0,1]$ and $\mathbb{R}$. For our purpose these measures
will always be the Lebesgue measure, which we will denote by
$\lambda$.
\begin{prop}$($Tietze's Extension Type$)$\label{Tidze_extension_obeying_null}
Let $E \subsetneq [0,1]$ be closed and let $f:E\mapsto [0,1]$ be a
continuous function that satisfy the property that for a measurable
set $A\subset E$, $\lambda(A)=0$ if and only if $\lambda(f(A))=0$.
Then there exists a continuous function $F:[0,1]\mapsto [0,1]$ such
that \\
$(i)$ $F_{\mid E}=f$,\\
$(ii)$ $F$ satisfies the property that for a measurable set
$A\subset [0,1]$, $\lambda(A)=0$ if and only if $\lambda(F(A))=0$.
\end{prop}
\begin{proof}
Since $E$ is closed it is a compact subset of $[0,1]$. Therefore $E$
has greatest and least members $m$ and $M$ respectively. If $m\neq
0$ or $M\neq 1$ then extend $f$ to a function $h$ on
$E_{1}=E\cup\{0\}\cup\{1\}$ by assigning the values $f(m)$ and
$f(M)$ at the points $0$ and $1$ respectively. The function $h$ is
continuous on $E_{1}$ and satisfies the same condition as $f$
relative to $E_{1}$. So without loss of generality we can assume
$0,1\in E$.\\
The complement of $E$ is a open set in $[0,1]$ and $E^{c}\subset
(0,1)$. Then $E^{c}$ can be written as a countable disjoint union of
intervals $\overset{\infty}{\underset{i=1}\cup}(a_{i},b_{i})$.
Then note that $a_{i},b_{i}\in E$ for all $i$.\\
So we only have to define an extension on $(a_{i},b_{i})$. Define
\begin{equation}
\nonumber
  F(x) = \begin{cases}
  f(x) &\text{ if }  x \in  E, \\
  \lambda f(a_{i})+(1-\lambda)f(b_{i}) &\text{ if } x=\lambda a_{i}+(1-\lambda)b_{i}\in (a_{i},b_{i}),\\
  &\indent 0 <\lambda <1 \text{ and }f(a_{i})\neq f(b_{i}),\\
  \frac{2(1-f(a_{i}))}{b_{i}-a_{i}}(x-a_{i})+f(a_{i}) &\text{ if
  }a_{i}<x\leq \frac{a_{i}+b_{i}}{2}\text{ and }\\
  &f(a_{i})=f(b_{i})<1,\\
  \frac{2(1-f(b_{i}))}{a_{i}-b_{i}}(x-b_{i})+f(b_{i}) &\text{ if
  }\frac{a_{i}+b_{i}}{2}\leq x< b_{i} \text{ and }\\
  &f(a_{i})=f(b_{i})<1,\\
  \frac{2(x-a_{i})}{a_{i}-b_{i}}+1 &\text{ if
  }a_{i}< x\leq \frac{a_{i}+b_{i}}{2} \text{ and }\\
  &f(a_{i})=f(b_{i})=1,\\
  \frac{2(x-b_{i})}{b_{i}-a_{i}}+1 &\text{ if
  }\frac{a_{i}+b_{i}}{2}\leq x<b_{i}  \text{ and }\\
  &f(a_{i})=f(b_{i})=1.\\
  \end{cases}
\end{equation}
The function $F$ is now continuous, as it is a linear interpolation
obtained from $f$ and the construction satisfy the required
conditions.
\end{proof}
\begin{thm}$($Foran,
\cite{MR807978}$)$\label{structure_theorem_for_Lusin_functions} A
necessary and sufficient condition for a continuous function
$F:[0,1]\mapsto [0,1]$ to satisfy condition $(N)$ relative to
$[0,1]$ is that there exists a sequence of measurable sets
$E_{n}\subseteq [0,1]$, $n=0,1,\cdots,$ such that the following properties are true:\\
$(i)$ $[0,1]=\overset{\infty}{\underset{n=0}\cup} E_{n}$,\\
$(ii)$ $\lambda(F(E_{n}))\leq n\lambda(E_{n})$ for all $n\geq 0$,\\
$(iii)$ for each $n>0$, $F$ is one to one on $E_{n}$.
\end{thm}
\begin{prop}\label{structure_of_measurable_functions}
Let $F:[0,1]\mapsto [0,1]$ be a measurable function such that for
any measurable set $A\subset [0,1]$, $\lambda(A)=0$ if and only if
$\lambda(F(A))=0$. Then there exists a measurable set $E\subseteq
[0,1]$ such that $\lambda(E)>0$ and $F$ is one to one on $E$.\\
Moreover, if $Y_{0}\subseteq [0,1]$ is such that $\lambda(Y_{0})>0$,
then there exists $Y_{1}\subseteq Y_{0}$ with $\lambda(Y_{1})>0$
such that $F$ is one to one on $Y_{1}$.
\end{prop}
\begin{proof}
Let $\epsilon >0$. By Lusin's theorem, choose a closed set $H\subset
[0,1]$ such that $\lambda([0,1]\setminus H)< \epsilon$ and $F_{\mid
H}$ is continuous relative to $H$. Clearly, $F_{\mid H}$ satisfy the
property that $A\subset H$, $\lambda(A)=0$ if and only if
$\lambda(F_{\mid H}(A))=0$. By Prop.
\ref{Tidze_extension_obeying_null}, extend $F$ to a continuous
function $\tilde{F}:[0,1]\mapsto [0,1]$ such that $\tilde{F}$ has
the property that for $A\subset [0,1]$, $\lambda(A)=0$ if and only
if $\lambda(\tilde{F}(A))=0$.\\
Now by Thm. \ref{structure_theorem_for_Lusin_functions}, choose
measurable subsets $E_{n}\subseteq [0,1]$ such that
$[0,1]=\overset{\infty}{\underset{n=0}\cup} E_{n}$,
$\lambda(\tilde{F}(E_{n}))\leq n\lambda(E_{n})$ for all
$n=0,1,\cdots$, and for each $n>0$, $\tilde{F}$ is one to one on
$E_{n}$.\\
Since $\lambda(\tilde{F}(E_{0}))=0$ so $\lambda(E_{0})=0$. If
$\lambda(E_{n}\cap H)=0$ for all $n>0$ then $\lambda(H)=0$, which is
not the case. Therefore there is a $n_{0}>0$ such that
$\lambda(E_{n_{0}}\cap H)>0$. But $\tilde{F}_{\mid E_{n_{0}}\cap
H}=F_{\mid E_{n_{0}}\cap H}$ and clearly $F$ is one to one on
$E_{n_{0}}\cap H$. Rename $E=E_{n_{0}}\cap H$.\\
Suppose $\lambda(Y_{0})>0$. By choosing $\epsilon>0$ small enough
one can make sure that the closed set $H$ in the first part of the
proof satisfies $\lambda(Y_{0}\cap H)>0$. The same argument as the
first part applies, and there exists a $n_{0}>0$ such that $F$ is
one to one on $Y_{1}=Y_{0}\cap H\cap E_{n_{0}}$ and
$\lambda(Y_{1})>0$.
\end{proof}

\bibliographystyle{plain}

 \end{document}